\documentclass[a4paper,11pt]{article}

\usepackage{theorem}
\usepackage[dvips]{graphics}
\newtheorem{definition}{Definition}
\newtheorem{proposition}{Lemma}
\newtheorem{theorem}{Theorem}
\newtheorem{cor}{Corollary}

\usepackage{graphics, amssymb, amsmath}
\def\disp{\displaystyle}

\def\C{\mathcal{C}}
\def\mb{\mathbb}
\def\mc{\mathcal}
 
\def\lan{\langle}
\def\ran{\rangle}

\begin{document}
\centerline{Utility maximization in a jump market model\footnote{Marie-Amelie Morlais\\ETH Zurich, Switzerland.\\
Ramistrasse 101, HG F27.3, 8006 Zurich,\\ 
Tel: +41 44 632 5859\\ e-mail: marieamelie.morlais@free.fr
}}

\begin{abstract}
In this paper, we consider the classical problem of utility maximization in a financial market allowing jumps. Assuming that the constraint set of all trading strategies
is a compact set, rather than a convex one, we use a dynamic method from which we derive a specific BSDE. To solve the financial problem, we first prove existence and uniqueness results for the introduced BSDE. This allows to give the expression of the value function and characterize optimal strategies for the problem.
\textbf{MSC classification (2000):} 91B28, 91B16, 60H10.
\end{abstract}

\section{Introduction and motivation}
The aim of this paper is to solve the utility maximization problem in a discontinuous framework: in this setting, the price process of stocks is assumed to be a local martingale driven by an independent one dimensional brownian motion and a Poisson point process.\\
\indent This problem is a classical one in the financial literature. But, as opposed to most of the papers dealing with the same problem (among them, we cite \cite{Fritelli} or \cite{Schach1}), we cannot rely on duality results, since we do not impose the convexity of the constraint set. The idea consists rather in adapting a dynamic method based on the dynamic programming principle: in the paper \cite{ImkelleretHu} dealing with the same topic, the authors already used this method to derive a specific BSDE. But, contrary to the present paper, they work in a brownian filtration and hence, results for quadratic BSDEs are available. Therefore, our main contribution is to establish new existence results for solution of the BSDE obtained in this discontinuous setting: the essential difficulty is to handle both the presence of jumps and the presence of a quadratic term. To this end and especially to handle the quadratic growth of the driver of the BSDE, we first apply the method introduced in the Brownian setting in \cite{mkobylanski}. The model is very analogous to \cite{Becherer} but, contrary to this last reference, we assume here that the price process has jumps and that there exists additional constraints on the portfolio. 
\indent The rest of the paper is structured as follows: in Section 2, we describe the market model by giving some preliminary remarks, specifying all the notations and natural assumptions, before introducing the utility maximization problem. Section 3 focuses on the theoretical results about the introduced BSDE. In Section 4, we go back to the financial problem and we apply the results established in Section 3. In a last section, we give an extension of the previous results by relaxing the assumption of compactness of the constraint set. Lengthy proofs are relegated to the last section.\\  

\section{The model and preliminaries}
\hspace{0.5cm} We consider a probability space ($\Omega$, $\mathbb{F}$, $\mathbb{P}$) equipped with two independent stochastic processes:
\begin{enumerate}
\item[.] A standard (one dimensional) brownian motion: $W$ =$(W_{t})_{t \in [0,T]}$.
\item[.] A real-valued
Poisson point process $p$ defined on $[0, T] \times \mb{R}\setminus{\{0\}}$. Referring to chapter 2 in \cite{Ikedawatanabe}, we denote by $N_{p}(ds,dx)$ the associated counting measure, such that its compensator is  $$\disp{\hat{N}_{p}(ds,dx) = n(dx)ds},$$
and the Levy measure $n(dx)$ (also denoted by $n$) is positive and satisfies
\begin{equation}\label{eq: conditionmeslevy}
\big( n(\{0\}) =0 \big) \;\textrm{and} \; \left( n((1 \wedge |x|)^{2}) := \disp{\int_{\mb{R}\setminus{\{0\}}}(1 \wedge |x|)^{2}n(dx) < \infty}\right). \end{equation}
 \end{enumerate}
 These two processes are considered on $[0, T]$: $T$ is called the horizon or maturity time in the financial context and, in all the sequel, it is fixed and deterministic.
For technical reasons, $n$ is assumed to be finite and,
in all the paper, we denote by $\mc{F}$ the filtration generated by the two processes $W$ and $N_{p}$ (and completed by $\mathcal{N}$, consisting in all the $\mb{P}$-null sets). 
Using the same notations as in \cite{Ikedawatanabe}, we denote by $\tilde{N}_{p}(ds,dx)$ ($\tilde{N}_{p}(ds, dx) := N_{p}(ds, dx) - \hat{N}_{p}(ds, dx)$) the compensated measure, which is a martingale random measure: in particular, for any predictable and locally square integrable process $F$,
the stochastic integral $F \cdot \tilde{N}_{p} := \disp{\int F_{s}(x) \tilde{N}_{p}(ds,dx)}$ is a locally square integrable martingale. \\
We denote by $Z \cdot W $ (resp. $U \cdot \tilde{N}_{p}$) the stochastic integral of $Z$ w.r.t. $W$ (resp. the stochastic integral of $U$ w.r.t. $\tilde{N}_{p}$). Furthermore, the filtration $\mc{F}$ has the predictable representation property: i.e., for any local martingale $K$ of $\mc{F}$, there exists two predictable processes $Z$ and $U$ such that  
$$\forall \; t, \quad K_{t} = K_{0} + \big(Z \cdot W \big)_{t} + \big( U \cdot \tilde{N}_{p}\big)_{t}. $$ 
(In Section 2.2, we provide a definition of the hilbertian spaces, in which these stochastic integrals are considered).\\
 \subsection{Description of the model}
The financial market consists in one risk-free asset (assumed to have zero interest rate) and one single risky asset, whose price process is denoted by $S$.
In particular, the stock price process is a one dimensional local martingale satisfying 
\begin{eqnarray}\label{eq: eqavecsauts} dS_{s} = S_{s-}\left( b_{s}ds + \sigma_{s}dW_{s} + \int_{\mathbb{R}^{*}}{\beta_{s}(x)\tilde{N}_{p}(ds,dx)} \right) .  \end{eqnarray}
All processes $b$, $\sigma$ and $\beta$ are assumed to be bounded and predictable and $\beta$ satisfies: $\beta > -1$. This last condition implies that the stochastic exponential $\mc{E}(\beta \cdot \tilde{N}_{p})$ is positive, $\mb{P}$-a.s.: hence, the process $S$ is itself almost surely positive.
The boundedness of $\beta$, $\sigma$ and $\theta$ ensures both existence and uniqueness results for the SDE (\ref{eq: eqavecsauts}). Then, provided that: $\sigma \neq 0$, we can define $\theta$ by: $\theta_{s} =\sigma_{s}^{-1}b_{s}$ ($\mathbb{P}$-a.s. and for all $s$).  
 The process $\theta$, also called market price of risk process, is supposed to be bounded and, under this assumption, the measure $\mb{P}^{\theta}$ with density
 $$ \dfrac{d\mathbb{P}^{\theta}}{d\mathbb{P}} = \mathcal{E}_{T}(-\int_{0}^{.}\theta_{s}dW_{s}), $$
 is a risk-neutral measure: this means that, under $\mb{P}^{\theta}$, the price process $S$ is a local martingale. \\
\indent In what follows, we introduce the usual notions of trading strategies and self financing portfolio, assuming that all trading strategies are constrained to take their values in a closed set denoted by $\mc{C}$. In a first step and to make easier the proofs, this set $\mc{C}$ is supposed to be compact (this assumption is relaxed in the last section). Due to the presence of constraints in this model with finite horizon $T$, not any $\mc{F}_{T}$-measurable random variable $B$ is attainable by using contrained strategies. In that context, we adress the problem of characterizing dynamically the value process associated to the exponential utility maximization problem (in the sequel, we denote by $U_{\alpha}$ the exponential utility function with parameter $\alpha$, which is defined on $\mb{R}$ by: $ U_{\alpha}(\cdot) = - \exp(- \alpha \cdot)$). 

\begin{definition}\label{tradingstrat}
A predictable $\mathbb{R}$-valued process $\pi$ is a self-financing trading strategy, if it takes its values in a constraint set $\mc{C}$ and if the process $X^{\pi, t, x}$ such that    \begin{equation}\label{eq: wealthproc}
\forall \; s \in [t,\;T], \quad X_{s}^{\pi, t, x}:= \disp{x + \int_{t}^{s}\pi_{s}\frac{dS_{s}}{S_{s-}}},
\end{equation}
 is in the space $\mc{H}^{2}$ of semimartingales (see chapter 4, \cite{Protter}). Such a process $X^{\pi}:=X^{\pi, t, x}$ stands for the wealth of an agent having strategy $\pi$ and wealth $x$ at time $t$.
 \end{definition}
\indent Now, as soon as the constraint set $\mc{C}$ is compact, the set consisting of all constrained strategies satisfies an additional integrability property.
\begin{proposition}\label{classequality}
 Under the assumption of compactness of the constraint set $\mathcal{C}$,
all trading strategies $\pi:=(\pi_{s})_{s \in [t, T]}$ as introduced in Definition \ref{tradingstrat} satisfy 
 \begin{equation}\label{eq: uniformint}   
\{\disp{ \exp(-\alpha X_{\tau}^{\pi}), \tau \; \mathcal{F}\textrm{-stopping time} \;  } \} \; \textrm{is a uniformly integrable family}.
 \end{equation}
\end{proposition}
We denote by $\mathcal{A}_{t}$ the admissibility set where the subscript $t$ indicates that we start the wealth dynamics at time $t$: more precisely, this set consists in all the strategies whose restriction to the interval $ [0, t] $ is equal to zero and which satisfy both Definition \ref{tradingstrat} and the condition (\ref{eq: uniformint}). In Section 4.3, we provide a justification of this last technical condition, which has already been used in \cite{ImkelleretHu} in a Brownian setting. The usual and much more restrictive admissibility condition consists in assuming that the wealth process $X^{\pi}$ is bounded from below (uniformly over all strategies $\pi$).\\
 \indent Before proving this result, we introduce another notion which can also be found in \cite{Dellacheriemeyer}: a martingale $M$ is said to be in the class of BMO martingales if there exists a constant $c$, $c > 0$, such that, for all $\mc{F}$-stopping time $\tau$,
$$   \textrm{ess}\displaystyle{\sup_{\Omega}\mathbb{E}^{\mathcal{F}_{\tau}}(\lan M \ran _{T}-\lan M \ran_{\tau} )} \leq c^{2}
\;\textrm{and} \; |\Delta M_{\tau}|^{2} \leq c^{2}.   $$
(In the continuous case, the BMO property follows from the first condition, whereas, in the discontinuous setting, we need to ensure the boundedness of the jumps of $M$).
 The following result, referred as Kazamaki's criterion and also stated in \cite{Kazamaki}, relates the martingale property of a stochastic exponential to a BMO condition.
\begin{proposition}\label{kamazaki} \textbf{(Kazamaki's criterion)} $\quad$Let $\delta$ be such that: 0 $< \delta < \infty $ and $M$ a BMO martingale satisfying: $\Delta M_{t} \geq -1 + \delta$, $\mathbb{P}$-a.s. and for all $t$, then $\mc{E}(M)$ is a true martingale.
\end{proposition}
\textbf{Proof of Lemma \ref{classequality}} \hspace{0.5cm} We first rely on the expression (\ref{eq: wealthproc}) of $X^{\pi}$ and on the definition of $\theta$ to claim
$$\disp{ dX_{t}^{\pi} = \pi_{t} \sigma_{t}\theta_{t}dt + \pi_{t} \sigma_{t}dW_{t} +  \int_{\mb{R}\setminus\{0\}}\pi_{t} \beta_{t}\tilde{N}_{p}(dt, dx)}. $$
Applying a generalized It\^o's formula to $U = e^{- \alpha X^{\pi}}$ (one reference is Theorem 5.1, Chapter 2 in \cite{Ikedawatanabe}), we get
\begin{tabbing}
$dU_{t} =$ \= $ U_{t}  \big(- \alpha\pi_{t} \sigma_{t}dW_{t} +  \disp{\int_{\mb{R}\setminus\{0\}}(e^{- \alpha\pi_{t} \beta_{t}} - 1)\tilde{N}_{p}(dt, dx) }\big)$\\
\\
 \> $+ U_{t}  \big(- \alpha \pi_{t} \sigma_{t}\theta_{t}dt + \frac{\alpha^{2}}{2}|\pi_{t} \sigma_{t}|^{2}dt + \disp{\int_{\mb{R}\setminus\{0\}}(e^{- \alpha\pi_{t} \beta_{t}} - 1 + \alpha\pi_{t} \beta_{t})n(dx)ds}). $\\
\end{tabbing}
 Using the definition of the stochastic exponential $\mc{E}(M)$ of $M$, $U$ has the multiplicative form
\begin{equation}\label{eq: expoalpha}
 U_{t} = U_{0}\mc{E}_{t}(\int_{0}^{.}- \alpha\pi_{s} \sigma_{s}dW_{s} + \int_{0}^{.}\int_{\mb{R}\setminus\{0\}}(e^{- \alpha\pi_{s} \beta_{s}} - 1) )\tilde{N}_{p}(ds, dx))e^{\bar{A}_{t}^{\pi}}, 
\end{equation}
where the process $\bar{A}^{\pi}$ is such that
$$ \bar{A}_{t}^{\pi} = \int_{0}^{t}\left( - \alpha \pi_{s} \sigma_{s}\theta_{s} + \frac{\alpha^{2}}{2}|\pi_{s} \sigma_{s}|^{2} + \int_{\mb{R}\setminus\{0\}}(e^{- \alpha\pi_{s} \beta_{s}} - 1 + \alpha\pi_{s} \beta_{s})n(dx)\right) ds.$$ 
Hence, $\bar{A}^{\pi}$ is a bounded process, thanks to the boundedness of the parameters of the SDE (\ref{eq: eqavecsauts}), the finiteness of the measure $n$ and the compactness of the constraint set $\mc{C}$).
Besides, thanks to Lemma \ref{kamazaki}, the stochastic exponential appearing in (\ref{eq: expoalpha}) is a true martingale and the desired property (\ref{eq: uniformint}) results from this decomposition.
\begin{flushright}
$\square$
\end{flushright}

 \subsection{Preliminaries}
In the sequel, we denote by $\mc{S}^{\infty}(\mb{R})$ the set of all adapted processes $Y$ with c\`adl\`ag paths (c\`adl\`ag stands for right continuous with left limits) such that: $\textrm{ess} \disp{\sup_{t, \omega}}|Y_{t}(\omega)| < \infty$, $L^{2}(W)$ denotes the set of all predictable processes $Z$ such that 
$$\disp{\mathbb{E}\left(\int_{0}^{T}|Z_{s}|^{2}ds \right) <\infty. } $$
  and $L^{2}(\tilde{N}_{p})$ denotes the set of all $\mathcal{P} \otimes \mathcal{B}(\mathbb{R}\setminus\{0\}) $-measurable processes $U$ such that 
$$ \disp{ \mathbb{E}\left(\int_{[0, T]\times \mathbb{R}^{*}}|U_{s}(x)|^{2}n(dx)ds \right) < \infty.}$$
$\mathcal{P}$ stands for the $\sigma$-field of all predictable sets of $[0, T] \times \Omega$ and $\mathcal{B}(\mathbb{R}\setminus\{0\})$ the Borel field of $\mb{R}\setminus\{0\}$.
We introduce $L^{0}(n)$ (denoted by $L^{0}(n, \mb{R}, \mb{R}\setminus\{0\})$ in \cite{Becherer}) as being the set of all the functions $u$, which map $\mb{R}$ in $\mb{R}\setminus\{0\}$. This set is equipped with the topology of convergence in measure. $L^{2}(n)$ (resp. $L^{\infty}(n)$) is the subset of all functions in $L^{0}(n)$ such that: $\disp{\mb{E}(\int_{0}^{T}|u(x)|^{2}n(dx))} < \infty$ (resp. $u$ takes bounded values).
  \\
\indent A solution of a BSDE with jumps, which is given by its terminal condition $B$ and its generator $f$, is a triple of processes ($Y$, $Z$, $U$) defined on
$ S^{\infty}(\mb{R}) \times L^{2}(W) \times L^{2}(\tilde{N}_{p} )$ such that: $\disp{\int_{0}^{T}|f(s, Y_{s}, Z_{s}, U_{s})|ds} < \infty, \; \mb{P}$-a.s. and satisfying 
\begin{equation}\label{eq: EDSRavecsauts} 
Y_{t} = B + \int_{t}^{T} f(s, Y_{s-}, Z_{s}, U_{s})ds - \int_{t}^{T} Z_{s}dW_{s} - \int_{t}^{T}\int_{\mathbb{R}^{*}}U_{s}(x)\tilde{N}_{p}(ds, dx).
\end{equation}
Throughout this paper, we study a specific BSDE with jumps whose terminal condition coincide with the contingent claim $B$, which is a bounded $\mathcal{F}_{T}$-measurable random variable, and whose generator $f$ is independent of $y$ (this is the case in the financial application).
Besides and since we do not work on a brownian filtration, the processes $Z$ and $U$ have to be predictable, for any solution of the BSDE (\ref{eq: EDSRavecsauts}) .\\
\indent A solution of such BSDEs with jumps is usually defined on
$ S^{2} \times L^{2}(W) \times L^{2}(\tilde{N}_{p} )$, $S^{2}$ being equipped with the following norm: $|Y|_{S^{2}} = \mathbb{E}(\disp{\sup_{t \in [0, T]}|Y_{t}|^{2}})^{\frac{1}{2}}$ (our main references for studies on BSDEs with jumps, are \cite{BarlesBuck} and \cite{Royer2}). But, in this paper, the results of the aforementionned papers cannot be applied, since the generator of the BSDE we are interested in, does not satisfy the usual conditions (it is not Lipschitz w.r.t. its variable $z$).

\section{The quadratic BSDE with jumps}
\subsection{Main assumptions}
In all the sequel (except in the proof of a priori estimates), we use the explicit form of the generator
\begin{eqnarray}\label{eq: generateurf}
 \disp{ f(s,z,u) = \displaystyle{\inf_{\pi \in
 \mathcal{C}}\left(\frac{\alpha}{2}|\pi\sigma_{s} - (z+ \frac{\theta_{s}}{\alpha})|^{2} + |u - \pi \beta_{s}|_ {\alpha}\right) -
  \theta_{s} z - \frac{|\theta_{s}|^{2}}{2\alpha}}},\\
\nonumber \end{eqnarray}
where the processes $\beta$, $\theta$ and $\sigma$ are defined in Section 2.1.
This expression of the generator will be justified in Section 4.3.  
We introduce the notation $|\cdot|_{\alpha}$ as being the convex functional such that 
\[ \begin{array}{ll}
\forall \; u \in (L^{2} \cap L^{\infty})(n),\; |u|_{\alpha} & =\disp{\int_{\mb{R}\setminus\{0\}}\frac{\exp(\alpha u(x)) - \alpha u(x) -1}{\alpha}n(dx)}, \\
\\
& =\disp{\int_{\mb{R}\setminus\{0\}}g_{\alpha}(u(x)) n(dx) },\\ 
\end{array}\] 
with the real function $g_{\alpha}$ defined by: $g_{\alpha}(y) = \frac{\exp(\alpha y) - \alpha y - 1 }{\alpha}$.
 As in \cite{Becherer}, we assume the finiteness of the measure $n$.
 In all that paper, $B$ is a bounded $\mc{F}_{T}$-measurable random variable and we use these two standing assumptions on the generator $f$
\begin{enumerate}
 \item[$(H_{1})$.] 
The first assumption denoted by ($H_{1}$) consists in specifying both a lower and an upper bound
\[  \begin{array}{l}
 \forall \; z, u \in \mb{R} \times (L^{2} \cap L^{\infty})(n) \\
\\
- \theta_{s} z - \frac{|\theta_{s}|^{2}}{2\alpha}  \le f(s, z, u) \le \frac{\alpha}{2}|z|^{2}  + |u|_{\alpha} ,\; \mb{P}\textrm{-a.s. and for all}\;s. \\
\end{array} \]
\item[($H_{2}$).] The second assumption, referred as ($H_{2}$), consists in two estimates: 
the first one deals with the increments of the generator $f$ w.r.t. $z$
\[  \begin{array}{l}
  \exists \; C > 0, \;\kappa \in BMO(W), \;\forall \; z, \; z' \in \mathbb{R},\; \forall u \in L^{2}(n(dx)),   \\
\\
|f(s, z, u) - f(s, z', u )| \leq  C(\kappa_{s} + |z| + |z'|)|z - z'| \\
\\
\end{array} \]
The second estimate deals with the increments w.r.t. $u$
\[  \begin{array}{l}
      \forall  z \in \mathbb{R}, \;\forall \;  u, u' \in (L^{2}\cap L^{\infty})(n(dx)),     \\
\\
f(s, z, u) - f(s, z, u') \leq \disp{\int_{\mb{R}\setminus{\{0\}}}\gamma_{s}(u, u')(u(x) - u'(x))n(dx)},\;  \\
    \end{array} \]
with the following expression for $\gamma_{s}(u, u^{'})$ for all $s$
\[ \begin{array}{l}
\disp{ \; \gamma_{s}(u, u')  \; =  }\\
\\
\displaystyle{\sup_{\pi \in \mathcal{C}} \left(\int_{0}^{1} g_{\alpha}^{'}(\lambda(u- \pi \beta_{s}) + (1 - \lambda)(u'- \pi\beta_{s})(x))d\lambda \right)}\mathbf{1}_{u \geq u'} \quad \quad \quad  \\
  \\
 \; \; + \;  \displaystyle{ \inf_{\pi \in \mathcal{C}} \left(\int_{0}^{1}g_{\alpha}^{'}(\lambda(u- \pi \beta_{s}) + (1 - \lambda)(u'- \pi \beta_{s})(x)d\lambda \right)}\mathbf{1}_{u < u'},\quad \quad \quad  \\
 \end{array}  \]
and for any fixed
 $s, \omega$,  this last expression makes sense as soon as $u$ and $u^{'}$ are in $(L^{2}\cap L^{\infty})(n(dx)) $. Now, for any arbitrary predictable processes $U$, $U^{'}$ taking their values in $ L^{2} \cap L^{\infty} (n)$), we can define a process\footnote{Since the increments of the generator $f$ are computed on the trajectories of predictable processes, (\ref{eq: processgamma}) is the expression of the process we are interested in, especially in the proof of the uniqueness result.} $\tilde{\gamma}$ by  

\begin{equation}\label{eq: processgamma}
\tilde{\gamma}_{s} =
\gamma_{s}(U_{s}, U^{'}_{s}), \end{equation}
and its predictability results from its expression in terms of both the predictable processes $U$, $ U^{'}$ and $\beta$.

\textbf{Comments on these assumptions}\\
$\bullet $ The expression (\ref{eq: generateurf}) of the generator $f$ clearly satisfies the assumption ($H_{1}$).
As in \cite{BriandetHu}, where the authors work in a brownian setting, this growth condition is useful to establish precise a priori estimates for any solution of the BSDE in $S^{\infty}(\mb{R}) \times L^{2}(W) \times L^{2}(\tilde{N}_{p})$. These estimates are given in lemma \ref{estim2} in the next paragraph and, analogously to \cite{ImkelleretHu}, their existence is the starting point of the proof of the main existence result.\\
\\
To check the first estimate in ($H_{2}$) we define, analogously as in \cite{ImkelleretHu} and for any $z, z^{'} \in \mb{R}$ and $u' \in (L^{2}\cap L^{\infty})(n(dx))$, 
 \[\left\{  \begin{array}{ll}
     \lambda_{s}(z, z^{'}) \;& :=  \frac{f(s, z, u) - f(s, z^{'}, u)}{z - z^{'}},  \; \textrm{if}\; z - z^{'} \neq 0,\\
     \\
    \lambda_{s}(z, z^{'}) &  := \; 0, \quad \textrm{otherwise}.\\
    \end{array} \right.
  \]  
It entails: $\lambda_{s}(z, z^{'}) \le C\big(\kappa_{s} + |z|+ |z^{'}| \big),$ 
and hence, the generator $f$ given by (\ref{eq: generateurf}) satisfies this first assumption with $C := \frac{\alpha}{2}$ and with $\kappa$ depending only on $\theta$ and on the two reals $\alpha$ and $ \disp{\sup_{\pi \in \mc{C}} |\pi|}$.\\
$\bullet $
To obtain the second estimate in ($H_{2}$), we rely on 
 \begin{tabbing}
 $ \quad$ \=$ \disp{f(s,z,u) -} \disp{ f(s,z,u') }$ \\
\\
 \> $ \disp{\leq \displaystyle{\sup_{\pi \in C}}\left( \int_{\mathbb{R}^{*}}g_{\alpha}((u - \pi\beta_{s})(x))- g_{\alpha}((u' -
 \pi\beta_{s})(x)n(dx)\right)}$ \\
\\
\> $ \disp{\leq \displaystyle{\sup_{\pi \in C}}\left( \int_{\mathbb{R}^{*}}\disp{\big(\int_{0}^{1}g_{\alpha}^{'}\big(\lambda(u - \pi\beta_{s}) + (1-\lambda)(u^{'} - \pi\beta_{s}) \big)  d\lambda}\big)(u -u^{'})(x) n(dx)\right)},$\\
\end{tabbing}
and we use that: 
$(u - u^{'})(x) = (u - u^{'})(x)\mathbf{1}_{u \geq u'} - (-(u -u^{'})(x))\mathbf{1}_{u < u^{'}}$,  \\
 and also:
$-\disp{\sup_{\pi }(-A^{\pi})} = \disp{\inf_{\pi }{A^{\pi}}},$ with $A^{\pi}$ such that
$$ A^{\pi}:= \disp{\left(\int_{0}^{1}g_{\alpha}^{'}\big(\lambda(u - \pi\beta_{s}) + (1-\lambda)(u^{'} - \pi\beta_{s}) \big)  d\lambda\right)}, $$
to claim that the expression given in ($H_{2}$) on page $7$ can be taken for $\gamma$.
To obtain the BMO property of the process given by (\ref{eq: processgamma}), we
use the compactness of $\mc{C}$ and we assume that both processes $U$ and $U'$ take their values $ L^{2} \cap L^{\infty}(n(dx))$ with $ |U_{s}|_{L^{\infty}(n)}, |U_{s}'|_{L^{\infty}(n)} \leq K$, to argue that
$$ \exists \; \delta_{K}, \bar{C}_{K} > 0,\; \textrm{s.t.} \quad  - 1 + \delta_{K} \leq \gamma_{s}(U_{s}, U_{s}^{'}) \leq \bar{C}_{K}, $$
which entails, in particular, that this process is in BMO($\tilde{N}_{p}$).
This BMO property is essential in the proof of the uniqueness result to justify the use of Girsanov's theorem.

\end{enumerate}
 \subsection{Theoretical results} 
\hspace{0.5cm} We now state the uniqueness and existence results for solutions of the BSDE (\ref{eq: EDSRavecsauts}) with parameters ($f$, $B$).
 \begin{theorem}\label{uniqueness} \textbf{(Uniqueness)} $\;$
$f$ given by (\ref{eq: generateurf}) and the terminal condition
$B$ being bounded, the BSDE of the form (\ref{eq: EDSRavecsauts}) with parameters ($f$, $B$) has at most one solution ($Y$, $Z$, $U$) in $ S^{\infty}(\mb{R}) \times L^{2}(W) \times L^{2}(\tilde{N}_{p})$.
\end{theorem}
\begin{theorem}\label{existence} \textbf{(Existence)} $\;$
$f$ given by (\ref{eq: generateurf}) and the terminal condition $B$ being bounded, there exists at least one solution of the BSDE (\ref{eq: EDSRavecsauts}) in $S^{\infty}(\mb{R}) \times L^{2}(W) \times L^{2}(\tilde{N}_{p})$. 
\end{theorem}

	     \subsection{A priori estimates }
 \hspace{0.5cm} In this part, we establish a priori estimates for any solution of BSDE (\ref{eq: EDSRavecsauts}) with a bounded terminal condition $B$ and with a generator $g$, which may be different from $f$ but satisfies ($H_{1}$). 

	\begin{proposition}\label{estim2}
For any BSDE of the form (\ref{eq: EDSRavecsauts}) with a generator $g$ satisfying ($H_{1}$) and a bounded terminal condition $B$, there exists three constants $C_{1}$, $C_{2}$ and $C_{3}$ depending only on $|B|_{\infty}$, $|\theta|_{S^{\infty}(\mb{R})}$ and $\alpha$, and such that, for any solution ($Y$, $Z$, $U$) in $S^{\infty}(\mb{R}) \times L^{2}(W) \times L^{2}(\tilde{N}_{p})$ and for any $\mathcal{F}$-stopping time $\tau$, $\tau$ taking its values in $[0, T]$,
	\[       \begin{array}{l}  
	    (i) \; \mb{P}\textrm{-a.s. and for all} \;\; t, \; t \in [0, T], \;\; C_{1} \leq
	     Y_{t} 
	     \leq C_{2},\;   \\
\\
	     (ii) \;   \mathbb{E}^{\mathcal{F}_{\tau}}\disp(\int_{\tau}^{T}|Z_{s}|^{2}ds + \int_{\tau}^{T}{\int_{\mathbb{R}^{*}}{|U_{s}(x)|^{2}
n(dx)}ds} \disp) \leq C_{3}.   \\

\\
	     \nonumber \end{array}  \]
	       \end{proposition}

\begin{cor}\label{equivalence}
Under the same assumptions than in lemma \ref{estim2} on the parameters $g$ and $B$ and for any solution ($Y$, $Z$, $U$) in $S^{\infty}(\mb{R}) \times L^{2}(W) \times L^{2}(\tilde{N}_{p})$ of the BSDE (\ref{eq: EDSRavecsauts}), 
\begin{itemize}
 \item there exists a predictable version $\tilde{U}$ of $U$ such that: $\tilde{U} \equiv U $ (in $L^{2}(\tilde{N}_{p})$). Noting $U$ instead of $\tilde{U}$, this process satisfies
$$ (i)(a) \; |U_{s}|_{L^{\infty}(n)} \leq 2 |Y|_{S^{\infty}(\mb{R})}, \;\textrm{and } \; (i)(b)\; \forall \; s, \;  |U_{s}|_{L^{2}(n)}^{2} \leq 4n(\mb{R}\setminus\{0\})|Y|_{S^{\infty}(\mb{R})}^{2}. $$
\item This equivalence result is satisfied
\begin{eqnarray}\label{eq: relationeq}
\nonumber \exists \; C> 0 \;  \frac{1}{C} \mathbb{E}\int_{[0, T] \times \mb{R}\setminus{\{0\}}}|U_{s}(x)|^{2}n(dx)ds  \leq \disp{\mathbb{E}\int_{0}^{T}|U_{s}|_{\alpha}ds} \\
\quad \quad \quad \quad \leq \disp{C \mathbb{E}\int_{[0, T] \times \mb{R}\setminus{\{0\}}}|U_{s}(x)|^{2}n(dx)ds,} \end{eqnarray}
with a constant $C$ depending only on $\alpha$ and $|Y|_{S^{\infty}(\mb{R})}$.
\end{itemize}

\end{cor} 

\paragraph*{Proof of the corollary}
To justify the first assertion (i)(a), we write
\begin{equation}\label{eq: jumpsofY}
|\Delta Y_{t}|= |Y_{t} - Y_{t-}| = \int_{\mb{R}\setminus{\{0\}}} |U_{t}(x)| N_{p}(\{t\}, dx),
\end{equation}
 and referring to the same notations as in chapter 2, \cite{Ikedawatanabe}, we have 
\begin{equation}\label{eq: formedessauts}\disp{\int_{\mb{R}\setminus\{0\}}|U_{t}(x)| N_{p}(\{ t\}, dx) = |U_{t}(p(t))|\mathbf{1}_{(t \in D_{p} )}},\end{equation}
 (the random set $D_{p}$ stands for the jump times of the Poisson point process).
 Defining $\tilde{U}$ by
$$\tilde{U}_{t}(x) = U_{t}(x)\mathbf{1}_{ |U_{t}(x)| \le 2 |Y|_{S^{\infty}(\mb{R})}}, $$ we obtain a predictable process, and it results from (\ref{eq: jumpsofY}) and (\ref{eq: formedessauts}) that:\\ $\tilde{U}_{t}(p(t))\mathbf{1}_{(t \in D_{p} )} = U_{t}(p(t))\mathbf{1}_{(t \in D_{p} )}$.
As a consequence
\[\begin{array}{l}
 \disp{\mb{E}\left(\int_{0}^{T}\int_{\mb{R}\setminus{\{0\}}} |(\tilde{U} - U)(s, x)|^{2} \hat{N}_{p}(ds, dx) \right)} \\
\\
\disp{ \quad \quad = \mb{E}\left(\sum_{s \in D_{p}, \; s \le T} |U(s, p(s)) - \tilde{U}(s, p(s))|^{2} \right)= 0,}\\
\end{array} \]
from which we deduce: $U = \tilde{U}$ (in $L^{2}(\tilde{N}_{p})$). We can now assume that, for any solution ($Y, Z, U$), $U$ satisfies: $ |U_{t}|_{L^{\infty}(n)} \leq 2 |Y|_{S^{\infty}(\mb{R})}, \; \mb{P}\textrm{-a.s. and for all} \; t.$
To justify the estimate in $L^{2}(n)$, we rely both on Cauchy Schwarz inequality and on the finiteness of the Levy measure. As soon as $U$ takes its values in $ L^{2} \cap L^{\infty}(n)$, the equivalence result (\ref{eq: relationeq}) holds. \\
\indent This last result is useful to prove Lemma \ref{convmonotone} also referred as ``monotone stability'' result and inspired by the work of Kobylanski in the Brownian setting in \cite{mkobylanski}.\\ 
\newline
\textbf{Proof of Lemma \ref{estim2}} 
 To establish (i), we assume the existence of a solution of the BSDE (\ref{eq: EDSRavecsauts}) with parameters ($g$, $B$) and we apply It\^o's formula to $e^{\alpha Y}$\\
    \begin{tabbing}
$ \disp{e^{\alpha Y_{t}} - e^{\alpha B} = }$ \= $\disp{\int_{t}^{T} \alpha e^{\alpha Y_{s}}\left( g(s,Z_{s},U_{s}) - \frac{\alpha}{2}|Z_{s}|^{2} - |U_{s}|_{\alpha}
\right)ds }$\\
\\
\> $\disp{ -\int_{t}^{T}\alpha e^{\alpha Y_{s}}Z_{s}dW_{s} - \int_{t}^{T}\int_{\mathbb{R}^{*}}e^{\alpha Y_{s}}(e^{\alpha U_{s}(x)} -1)\tilde{N}_{p}(dx,ds) }.$ \\
\end{tabbing}
 However, to be more rigorous, we should apply the standard procedure of localization, i.e we should first introduce the sequence of stopping times
$$ \tau^{m} = \disp{\inf_{0 \leq t \leq T} \big(\{ \int_{0}^{t}e^{2 \alpha Y_{s}}|Z_{s}|^{2}ds \geq m \} \cup \{\int_{0}^{t}\int_{\mb{R}\setminus\{0\}}e^{2 \alpha Y_{s}}|e^{ \alpha U_{s}(x)} -1|^{2}n(dx)ds \geq m \}\big)}, $$
and then take the conditional expectation w.r.t. $\mc{F}_{t}$, before passing to the limit as $m$ goes to $\infty$.
Hence and without loss of generality, we just take the conditional expectation w.r.t. $\mathcal{F}_{t}$ in It\^o formula applied to $e^{\alpha Y}$ and use the upper bound in ($H_{1}$) to obtain\\
$$\disp{ e^{\alpha Y_{t}} -  \mathbb{E}^{\mathcal{F}_{t}}(e^{\alpha B}) 
 \leq \mathbb{E}^{\mathcal{F}_{t}}(\int_{t}^{T} \alpha e^{\alpha Y_{s}}\underbrace{( g(s,Z_{s},U_{s}) -\frac{\alpha}{2}|Z_{s}|^{2} - |U_{s}|_{\alpha})}_{\leq 0}ds) }.$$
Since: $\mathbb{E}^{\mathcal{F}_{t}}(e^{\alpha B}) \leq 
e^{\alpha|B|_{\infty}}$, the right-hand side in (i) holds true
 with $C_{2} = |B|_{\infty}$. To obtain the left-hand side in (i), we introduce $A$ such that $A:= (g(s, Z_{s}) + ( Z_{s} \theta_{s} + \frac{|\theta_{s}|^{2}}{2 \alpha}))$: such a process $A$ is non negative, thanks to ($H_{1}$). Hence, we have   
 \begin{tabbing}
$Y_{t}  $\= $ \; = \; \disp{ Y_{T} + \int_{t}^{T} g(s, Z_{s}, U_{s})ds - \int_{t}^{T}Z_{s}dW_{s}-  \int_{t}^{T}\int_{\mathbb{R}^{*}}U_{s}(x)\tilde{N}_{p}(ds, dx) }$\\
 \> $ $\\
\> $ \;=  \;\disp{ Y_{T} +\; (A_{T} - A_{t}) - \int_{t}^{T} ( Z_{s} \theta_{s} + \frac{|\theta_{s}|^{2}}{2 \alpha})ds}$\\
\> $ $\\
 \> $ \quad \quad \; \disp{- \int_{t}^{T}Z_{s}dW_{s}-  \int_{t}^{T}\int_{\mathbb{R}^{*}}U_{s}(x)\tilde{N}_{p}(ds, dx)}.$\\
	   \end{tabbing}
We now introduce the density of a new probability measure $\mathbb{P}^{\theta}$ 
$$ d\mathbb{P}^{\theta} = \mathcal{E}_{T}(-\int_{0}^{.} \theta_{s}dW_{s})d\mathbb{P}.$$
 The BMO property of $- \theta \cdot W$ entails that the continuous stochastic exponential $\mathcal{E}(-\int_{0}^{.} \theta_{s}dW_{s})$ is a true martingale. Hence, $\mathbb{P}^{\theta}$ is an equivalent probability measure under which the process $W^{\theta} := W + \int_{0}^{.}\theta_{u}du$, is a standard Brownian motion. Finally, we get
\begin{equation}\label{eq: controletaun}
\disp{ Y_{t} = Y_{T} + (A_{T} -A_{t})- \int_{t}^{T}\frac{|\theta_{s}|^{2}}{2 \alpha} ds - \int_{t}^{T}Z_{s}dW_{s}^{\theta} - \int_{t}^{T}\int_{\mathbb{R}^{*}}U_{s}(x)\tilde{N}_{p}(ds, dx)}. \end{equation} Under $\mathbb{P}^{\theta}$, the two last terms are local martingales bounded from below. Using a standard localization procedure, we obtain the existence of ($\tau_{n}$) such that the two local martingales stopped at time $\tau_{n} $ are $\mc{F}$-martingales.
 Taking the conditional expectation w.r.t. $\mathcal{F}_{t}$ and under $\mb{P}^{\theta}$, we get
$$ Y_{t \wedge \tau_{n}} \ge \mb{E}^{\theta}\left( Y_{T \wedge \tau_{n}}- \int_{t}^{T}\frac{|\theta_{s}|^{2}}{2 \alpha} ds | \mc{F}_{t}\right),  $$
 $\mb{E}^{\theta}$ standing for the expectation under $\mb{P}^{\theta}$. Applying now the bounded convergence theorem and relying on the boundedness assumption on $\theta$ (given by ($H_{1}$)), we can let $n$ tend to $ +\infty$ in (\ref{eq: controletaun}) to claim
\[ Y_{t} = \mb{E}^{\theta}\big(Y_{t}|\mc{F}_{t}\big) \geq -|B|_{\infty} - \frac{|\theta|_{S^{\infty}}^{2}T}{2 \alpha} , \; \mb{P}^{\theta}\textrm{-a.s.} \]
This holds true $\mathbb{P}$-a.s., because of the equivalence of $\mathbb{P}^{\theta}$ and $\mathbb{P}$ and finally, the assertion
(i) is satisfied with $C_{1}$ such that	   
	  $$ C_{1} = -|B|_{\infty} - \frac{|\theta|_{S^{\infty}}^{2}T}{2 \alpha}.$$ 
\newline
\indent To prove assertion (ii), we apply It\^o's formula to ($Y - C_{2}$)$^{2}$ ($C_{2} := |B|_{\infty}$ is the upper bound in (i)). Since $Z \cdot W$ and $U \cdot \tilde{N}_{p}$ are true martingales (as square integrable stochastic integrals), we next integrate this formula between an arbitrary stopping time $\tau$ and $T$ and take the conditional expectation w.r.t. $\mathcal{F}_{\tau}$ to get
   \begin{tabbing}
 $\mathbb{E}^{\mathcal{F}_{\tau}}((Y_{\tau} - C_{2})^{2} $ \= $ - (Y_{T} - C_{2})^{2}) = \mathbb{E}^{\mathcal{F}_{\tau}}\big(\disp{\int_{\tau}^{T}{2(Y_{s}-C_{2})g(s,Z_{s},U_{s})ds}}\big) $\\
  \\
   \> $\quad - \mathbb{E}^{\mathcal{F}_{\tau}}\left(\disp{\int_{\tau}^{T}|Z_{s}|^{2}ds} + \disp{\int_{\tau}^{T}{\int_{\mathbb{R}^{*}}{|U_{s}(x)|^{2} n(dx)}ds}}\right).$ \\
  \end{tabbing}
To give an upper bound of $\big(Y_{s}-C_{2})g(s,Z_{s},U_{s}\big)$, we first rely on 
\[ \textrm{and}\; \left\{ \; \begin{array}{l} 
   0 \ge (Y - C_{2}) \ge -\big( 2|B|_{\infty} + \frac{|\theta|_{S^{\infty}}^{2} T}{2 \alpha}\big), \;   \\
\\
g(s,Z_{s},U_{s}) \ge  - \theta_{s} Z_{s} - \frac{|\theta_{s}|^{2}}{2 \alpha},\\
   \end{array} \right.  \]
to deduce the existence of a constant $C$ depending only on $\theta$ and $\alpha$ such that
\begin{equation}\label{eq: Controlestim2}
 (Y_{s}-C_{2})g(s,Z_{s},U_{s}) \le C\big(|\theta_{s} Z_{s}| + \frac{|\theta_{s}|^{2}}{2 \alpha} \big).
\end{equation}
 We now split the discussion into two cases:\\
$\bullet \;$if $\theta \equiv 0$, then $g$ is bounded by 0 from below and it implies 
\begin{tabbing}
 $\mathbb{E}^{\mathcal{F}_{\tau}}((Y_{\tau} - C_{2})^{2} $ \=$ - (Y_{T} - C_{2})^{2}) \le $\\
   \> $- \mathbb{E}^{\mathcal{F}_{\tau}}\left(\disp{\int_{\tau}^{T}|Z_{s}|^{2}ds} + \disp{\int_{\tau}^{T}{\int_{\mathbb{R}^{*}}{|U_{s}(x)|^{2} n(dx)}ds}}\right),$ \\
  \end{tabbing}
which leads to
$$\mathbb{E}^{\mathcal{F}_{\tau}}\left(\disp{\int_{\tau}^{T}|Z_{s}|^{2}ds} + \disp{\int_{\tau}^{T}{\int_{\mathbb{R}^{*}}{|U_{s}|^{2} n(dx)}ds}}\right) \le |Y_{T} - C_{2}|^{2} \le 4 |Y|_{S^{\infty}}.$$
$\bullet \;$ else if $\theta \neq 0$, we use: $ab \leq \frac{1}{2}(a^{2} + b^{2})$ in (\ref{eq: Controlestim2}) to argue the existence of a bounded process $C_{1}(\cdot)$ such that
\[2(Y_{s} - C_{2})g(s,Z_{s},U_{s})  \leq  C_{1}(s) +
\frac{1}{2}|Z_{s}|^{2}.  \]
From straightforward computations, we get 
\[ \begin{array}{l}
 \frac{1}{2}\left(\mathbb{E}^{\mathcal{F}_{\tau}}\left(\disp{\int_{\tau}^{T}|Z_{s}|^{2}ds} + \disp{\int_{\tau}^{T}{\int_{\mathbb{R}^{*}}{|U_{s}|^{2} n(dx)}ds}}\right) \right) \\
\\
\quad \quad \quad \quad \le \disp{\mathbb{E}^{\mathcal{F}_{\tau}}\left( \int_{\tau}^{T} C_{1}(s) ds + (Y_{T} - C_{2})^{2}  \right).}  \end{array} \]
and hence, (ii) follows.
  \begin{flushright}
$\square$
           \end{flushright} 
    \subsection{Uniqueness }
\textbf{Proof of Theorem \ref{uniqueness}} $\quad$
The idea consists in using a linearization procedure and in justifying, as in \cite{ImkelleretHu}, the application of Girsanov's theorem. We first consider ($Y^{1}$, $Z^{1}$, $U^{1}$) and ($Y^{2}$, $Z^{2}$, $U^{2}$) two solutions of the BSDE with jumps given by ($f$, $B$) and $C$ any positive constant such that: $|Y^{i}|_{S^{\infty}(\mb{R})} \leq C$ (thanks to (i) in Lemma \ref{estim2}, we can set $C:= |C_{1}| + |C_{2}|$). Then, we introduce $\hat{Y}$, $ \hat{Z}$ and $\hat{U} $  \\
\[ \hat{Y} = Y^{1} - Y^{2}, \; \hat{Z} =Z^{1} - Z^{2}, \;  \hat{U} = U^{1} - U^{2}.\]
 $\tau$ being an arbitrary $\mc{F}$-stopping time, we apply It\^o's formula to $\hat{Y}$ between $t \wedge \tau$ and $\tau$,
\[ \begin{array}{ll}
 \hat{Y}_{t \wedge \tau} - \hat{Y}_{\tau} = & \disp{\int_{t \wedge \tau }^{\tau}(f(s,Z_{s}^{1},U_{s}^{1}) - f(s,Z_{s}^{2},U_{s}^{2}))ds }\\
 \\
 & \; \; \disp{-\int_{t \wedge \tau}^{\tau}\hat{Z}_{s}dW_{s} -\int_{t \wedge \tau}^{\tau}\int_{\mathbb{R}^{*}}{\hat{U}_{s}\tilde{N}_{p}(ds,dx)}. }  \end{array}\]
Since the generator does not satisfy the usual conditions: it is not Lipschitz neither in $z$ nor in
$u$), we rely on assumption ($H_{2}$), i.e. on the estimates of the increments of the generator $f$ 
We first introduce $\lambda:= \lambda(Z^{1},Z^{2})$ 
\[\left\{  \begin{array}{ll}
     \lambda_{s}(Z_{s}^{1}, Z_{s}^{2}) \;& =  \frac{f(s,Z_{s}^{1},U_{s}^{1}) - f(s,Z_{s}^{2},U_{s}^{1})}{Z^{1}_{s} - Z_{s}^{2}},  \; \textrm{if}\; \hat{Z} \neq 0,\\
     \\
    \lambda_{s}(Z_{s}^{1}, Z_{s}^{2}) &  = \; 0, \quad \textrm{else}.\\
    \end{array} \right.
  \]  
Thanks to ($H_{2}$), we get
\begin{tabbing}
$  f(s,Z_{s}^{1}  $\= $,U_{s}^{1}) - f(s,Z_{s}^{2},U_{s}^{2})$\\
$ $\\
\> $ = f(s,Z_{s}^{1},U_{s}^{1}) - f(s,Z_{s}^{2},U_{s}^{1}) + f(s,Z_{s}^{2},U_{s}^{1}) - f(s,Z_{s}^{2},U_{s}^{2}) $\\
\> $   $ \\
\> $ \leq \lambda_{s}(Z_{s}^{1},Z_{s}^{2} ) \hat{Z}_{s} + \disp{\int_{\mb{R^{*}}}\gamma_{s}( U_{s}^{1}, U_{s}^{2})\hat{U}_{s}(x)n(dx)},$ \\
\end{tabbing}
where the process $\tilde{\gamma} = (\gamma_{s}( U_{s}^{1}, U_{s}^{2}))$ is defined analagously as in (\ref{eq: processgamma}) and for all $s$ by 
\[ \begin{array}{ll}
\disp{ \; \tilde{\gamma}_{s}  \; = } & 

\displaystyle{\sup_{\pi \in \mathcal{C}} \left(\int_{0}^{1} g_{\alpha}^{'}(\lambda(U_{s}^{1}- \pi \beta_{s}) + (1 - \lambda)(U_{s}^{2}- \pi\beta_{s})(x))d\lambda \right)}\mathbf{1}_{u \geq u'} \quad \quad \quad  \\
  \\
 \; \; & +\;  \displaystyle{ \inf_{\pi \in \mathcal{C}} \left(\int_{0}^{1}g_{\alpha}^{'}(\lambda(U_{s}^{1}- \pi \beta_{s}) + (1 - \lambda)(U_{s}^{2}- \pi \beta_{s})(x)d\lambda \right)}\mathbf{1}_{u < u'}.\quad \quad \quad  \\
 \end{array}  \]
As in subsection 3.1, $\tilde{\gamma}$ is explicitely expressed in terms of both the two predictable processes $U^{i}$ ($i=1, 2$) and the predictable process $\beta$ and hence,
it is itself predictable. Thanks to Corollary \ref{equivalence}, we also obtain that for $i:$=1, 2, $|U_{s}^{i}|_{L^{\infty}(n)} \leq 2|Y^{i}|_{S^{\infty}(\mb{R})} \leq 2 C,$ ($C$ is given at the beginning of this proof), which implies 
$$\exists \; \delta > 0, \; \bar{C} > 0,\; \textrm{s.t.} \; -1 + \delta \leq \gamma_{s}(U_{s}^{1}, U_{s}^{2}) \leq \bar{C}, \; \mb{P}\textrm{-a.s.} \; \textrm{and for all}\;  s.$$
Since the process $\lambda$ satisfies
$$|\lambda_{s}(Z_{s}^{1}, Z_{s}^{2})| \leq C(\kappa_{s} + |Z_{s}^{1}| + |Z_{s}^{2}|),$$
the BMO properties of $\lambda \cdot W$ and $\tilde{\gamma} \cdot \tilde{N}_{p}$ result from the assumption $(H_{2})$ and also from the BMO property of $\int_{0}^{.} Z^{i}dW_{s}$, which holds true for $i$ = 1, 2 (we refer to (ii) in Lemma \ref{estim2}).
Defining $M^{1}$ and $M^{2}$ such that
$$M^{1} = \lambda \cdot W \quad \textrm{and }\quad  M^{2}  = \tilde{\gamma} \cdot \tilde{N}_{p},$$
and setting: $d\mathbb{Q}$ = $\mathcal{E}(M^{1} + M^{2})d\mathbb{P}$, it results from Kazamaki's criterion that $\mc{E}(M^{1} +M^{2})$ is a true martingale and hence, $\mathbb{Q}$ is an equivalent probability measure. Using Girsanov's theorem, we set: $W^{\lambda}:= W - \lan{W}, \lambda \cdot W \ran$ and $\tilde{N}^{\gamma} := \tilde{N}_{p} - \lan{\tilde{N}_{p}, \tilde{\gamma} \cdot \tilde{N}_{p}} \ran $ and we introduce  
$$\forall \; t, \quad M_{t} = \int_{0}^{t}\hat{Z}_{s}dW^{\lambda}_{s} + \int_{0}^{t}\int_{\mathbb{R}^{*}}\hat{U}_{s}(x) \tilde{N}^{\gamma}(ds,dx),  $$ 
which is a local martingale under $\mb{Q}$. We conclude by relying on a standard localization procedure: there exists a sequence $(\tau^{n})$ of stopping times converging by increasing to $T$ and such that: $M_{t\wedge \tau^{n}}$ is a martingale. Hence, under the measure $\mb{Q}$, $ \hat{Y}_{\cdot \wedge \tau^{n}}$ is a submartingale, which yields
\[ \hat{Y}_{t \wedge \tau^{n}} \le \mb{E}^{\mb{Q}}\big(\hat{Y}_{\tau^{n}} |\mc{F}_{t}\big). \] Using the bounded convergence theorem, we finally get: $\hat{Y}_{t} \leq 0,$ $\mathbb{Q}$-a.s. and $\mathbb{P}$-a.s., because of the equivalence of $\mb{P}$ and $\mb{Q}$. Thanks to the symmetry of this problem, we can conclude: $\hat{Y} = 0$.
\begin{flushright}
$\square$
           \end{flushright} 
    \subsection{Existence }
\hspace{0.5cm} To prove the existence for a solution of the BSDE (\ref{eq: EDSRavecsauts}) with parameters ($f$, $B$) and with $f$ given by (\ref{eq: generateurf}), we proceed in three main steps and, for more convenience, we provide here a description of each step:\\
$\bullet \;$ The first step consists in the construction of an approximating sequence $(f^{n})$ of $f$ such that each generator is lipschitz. \\
$ \bullet \;$ In a second step and referring to already known results, we justify the existence of solutions to the BSDEs given by ($f^{n}$, $B$) and we provide precise estimates.\\
$ \bullet \;$
The last step consists in justifying a stability result for the solutions of these BSDEs (analogous to the one provided by \cite{mkobylanski}) and deduce from it the existence of a limit, which solves the original BSDE.\\ 
The proof being constructive, we need the explicit form (\ref{eq: generateurf}) of the generator,
 satisfying in particular the assumptions ($H_{1}$) and ($H_{2}$) stated in section 3.1.

\subsubsection{Step 1: Approximation by a truncation argument}
\hspace{0.5cm} 
Our aim is to construct an increasing sequence ($f^{m}$), such that each $f^{m}$ satisfies the assumption ($H_{1}$) with the same parameters as $f$. For all $m$, $m \geq M$, with $M$ given by $M:= 2(|C_{1}| + |C_{2}|)$ ($C_{1}$ and $C_{2}$ are the constants given in (i) by Lemma \ref{estim2}),
we define $f^{m}$ 
\[
\begin{array}{l}
 f^{m}(s, z, u) =  \\
\quad \disp{\inf_{\pi \in \mathcal{C}}\left( \frac{\alpha}{2}|\pi \sigma_{s} - (z + \frac{\theta_{s}}{\alpha})|^{2}\rho_{m}(z) + \int_{\mathbb{R}^{*}}g_{\alpha}(u - \pi \beta_{s})(x)\rho_{M }(u)(x)n(dx)\right)}  \\
\\
\quad  \quad  - \theta_{s} z - \frac{|\theta_{s}|^{2}}{2\alpha}. \\
  \end{array} \]
 The truncation function $\rho_{m}$, assumed to be continuously differentiable, is such that: \\
\textbullet $\;$ The sequence ($\rho_{m}$) is increasing with respect to $m$.\\
\textbullet $\;$ $\rho_{m}(x):$ = 1, if $|x| \leq m$, $0 \le \rho_{m}(x) \le 1 $, if $m \le |x| \le m + 1$, and  $\rho_{m}(x)$ = 0, if $|x| \geq m+1$. \\
\indent Now, we list and check the properties of $f^{m}$, required first to establish existence and uniqueness results for the BSDEs given by ($f^{m}$, $B$)
 and also, to ensure the passage to the limit, as $m$ goes to $\infty$.    \\
\begin{enumerate}
\item Each generator $f^{m}$ has the Lipschitz property: i.e, for any $m \geq M$,
\[   \begin{array}{l}
 \disp{\exists \; C(m) > 0,\; \; \forall \;t, \;\forall \;z, z^{'} \in  \mb{R},\; \forall \; u, u^{'} \in L^{2} \cap L^{\infty}(n), }\\
\\
 \disp{  |f^{m}(t,z,u) - f^{m}(t,z^{'},u^{'})| \leq } \\
\\
 \disp{\quad \quad  C_{m}\left(|z-z^{'}| + (\int_{\mathbb{R}^{*}}(u(x) -
  u^{'}(x))^{2} n(dx) )^{\frac{1}{2}} \right). }\\
   \end{array} \]
To handle the increments w.r.t $z$ and $u$, we proceed by linearization analogously as in Section 3.1. The Lipschitz property follows from the boundedness of the increments, which itself follows from the truncation.

\item For each $f^{m}$ and using Theorem 2.5 in \cite{Royer2}, we prove that the monotonicity assumption ($H_{\textrm{comp}}$) and the associate condition ($A_{\gamma}$) holds true. Both these conditions result from the assumption ($H_{2}$) on the increments in $u$. In particular, for each $f^{m}$, ($H_{2}$) is satisfied with: $\gamma^{m} :=  \gamma$, for all $m$, and, for any $u$ in $L^{2}(n)$, $\rho_{M}(u)$ is in $L^{2} \cap L^{\infty}(n)$.
Hence, a comparison result holds true for the BSDEs given by ($f^{m}$, $B$).\\

\item We have  \begin{equation}\label{eq: controlunif}
                     \mathop{\sup_{m}|f^{m}(s, 0, 0)| } \leq \frac{|\theta_{s}|^{2}}{2\alpha},                
                                                                  \end{equation}
 which entails that $\disp{\mathop{\sup_{m}|f^{m}(s, 0, 0)| } }$ is in $L^{1}(ds \otimes d\mb{P})$.

\item The sequence $(f^{m})$ converges to $f$ in the following sense \\
$ \forall \;s, \; s\in [0, T], \;\forall \;z \in \mb{R}, \forall \;u \in L^{2} \cap L^{\infty}(n) $ , 
$$ \quad \big( f^{m}(s, z, u) \to  f(s, z, u), \; \mb{P}\textrm{-a.s.}\big) \;\; \textrm{as}\; m \to \infty. $$
 This convergence results from the truncation of the functionals of $z$ and $u$ defining $f^{m}$. Besides and thanks to the increasing property of ($\rho_{m}$) and the positiveness of the square functional involving $z$, ($f^{m}$) is itself increasing.\\
\end{enumerate}

\subsubsection{Step 2: Useful properties of this approximation}
\hspace{0.5cm} Referring to the results in \cite{Royer2} or in \cite{Pardoux}, we obtain existence and uniqueness for a solution ($Y^{m}$, $Z^{m}$, $U^{m}$) in $S^{2} \times L^{2}(W) \times L^{2}(\tilde{N}_{p})$ of the BSDEs given by ($f^{m}, \;B$). 
To give estimates of $Y^{m}$ in $S^{\infty} $, which may depend on $m$, we state an auxiliary result.
\begin{proposition}\label{estimcaslip}
 Let ($Y^{n}$, $Z^{n}$, $U^{n}$) be a solution in $S^{2} \times L^{2}(W) \times L^{2}(\tilde{N}_{p})$ of a BSDE of type ($\textrm{Eq2}$) given by ($g^{n}$, $\bar{B}$), with a generator $g^{n}$ $L_{n}$-Lipschitz and a terminal condition $\bar{B}$ bounded, we have 
\begin{equation}\label{eq: controlbornee}
\exists \; K(L_{n}, T) > 0, \;\forall \;t, \; |Y_{t}^{n}|^{2} \leq K(L_{n}, T)\mb{E}\left(|\bar{B}|^{2} + (\int_{t}^{T}|g^{n}(s, 0, 0)|dC_{s})^{2}| \mc{F}_{t} \right). \end{equation}
\end{proposition}
(the detailed proof, adapted from \cite{BriandetCoquet}, is given in Appendix A).\\
\indent The solution being in $S^{\infty}(\mb{R}) \times L^{2}(W) \times L^{2}(\tilde{N}_{p}) $, we would like to get free from the dependence in $m$ of the estimates of $Y^{m}$ in $S^{\infty}(\mb{R})$: to this end, we justify that the estimates in Lemma \ref{estim2} hold true for the solution ($Y^{m}$, $Z^{m}$, $U^{m}$) and for all $m$: in fact, each $f^{m}$ satisfies the assumption ($H_{1}$) with the same parameters as $f$ and hence, 
\begin{itemize}
\item  ($Y^{m}$) is uniformly bound in $S^{\infty}(\mb{R})$.

\item $(Z^{m})$ and $(U^{m})$ are uniformly bounded in their respective Hilbert spaces, i.e. $ L^{2}(W)$ and $L^{2}(\tilde{N}_{p})$.\\
\end{itemize}
\indent As in Corollary \ref{equivalence}, we obtain an equivalence result for ($U^{m}$)$_{m}$: more precisely, there exists a constant $C$ independent of $m$ and satisfying
\begin{eqnarray*} \frac{1}{C} \int_{[0, T] \times \mathbb{R}^{*}}|U_{s}^{m}(x)|^{2}n(dx)ds  \leq \disp{\int_{0}^{T}|U_{s}^{m}|_{\alpha}ds \leq C\int_{[0, T] \times \mathbb{R}^{*}}|U_{s}^{m}(x)|^{2}n(dx)ds . } \end{eqnarray*}

\indent We now justify the existence of processes $\tilde{Y}$, $\tilde{Z}$ and $\tilde{U}$, which are the respective limits in a specific sense of ($Y^{m}$)$_{m}$, ($Z^{m}$)$_{m}$ and ($U^{m}$)$_{m}$. Using the comparison result for BSDE with jumps given by \cite{Royer2} and justified in Step 1, ($Y^{m}$)$_{m}$ is increasing. Hence, we can define $\tilde{Y}$ 
$$\tilde{Y}_{s} = \displaystyle{\lim_{m} \nearrow (Y_{s}^{m})},   \quad \mathbb{P}\textrm{-a.s. and for all} \; s.  $$
 Since
($Z^{m}$) and ($U^{m}$) are uniformly bounded in their respective BMO spaces and, in particular, in $L^{2}(W)$ and $L^{2}(\tilde{N_{p}})$, we can extract from both sequences converging subsequences in the weak sense. We denote by $\tilde{Z}$ and  $\tilde{U}$ their respective weak limits.\\
\subsubsection{Step 3: Convergence of the approximation}
\hspace{0.5cm}In this last step, we prove the convergence of ($Y^{m}$, $Z^{m}$, $U^{m}$) to a solution of the BSDE (\ref{eq: EDSRavecsauts}) with parameters $f$ and $B$. To this end, we justify the passage to the limit as $m$ goes to $\infty$ in 
\begin{equation}\label{eq: eqapprox} Y_{t}^{m} = Y_{T}^{m} + \int_{t}^{T} f^{m}(s,Z_{s}^{m},U_{s}^{m})ds - \int_{t}^{T}Z_{s}^{m}dW_{s} -
  \int_{t}^{T}\int_{\mathbb{R}^{*}}U_{s}^{m}\tilde{N}_{p}(ds, dx).  \end{equation}
To achieve this, the essential step is to establish a ``monotone stability'' result, which is an adaptation of Proposition 2.4 in \cite{mkobylanski}.\\

 \begin{proposition}\label{convmonotone} 
Assuming that\\
(1) There exists a sequence ($f^{m})_{m}$ such that,
 for all $s$ and for all converging sequences $(z^{m})_{m}$ and $(u^{m})_{m}$ taking their values in $\mathbb{R}$ and $L^{2}( n(dx))$, with ($u^{m}$) uniformly bounded in $L^{\infty}(n)$, 
$$\; \; \; \disp{\lim \nearrow f^{m}(s, z^{m}, u^{m}) :=  f(s, z, u)}, \quad \textrm{as} \; m \to \infty. $$
(2) Each $f^{m}$ satisfies $(H_{1})$ with the same parameters as $f$.\\
(3) There exists ($Y^{m}$, $Z^{m}$, $U^{m}$)$_{m}$ in $S^{\infty}(\mb{R}) \times L^{2}(W) \times L^{2}(\tilde{N}_{p})$, which are solutions of the BSDEs given by ($f^{m}$, $B$). \\
Then, ($Y^{m}$, $Z^{m}$, $U^{m}$) converges to ($\tilde{Y}, \tilde{Z}, \tilde{U}$) in the following sense
\begin{equation}\label{eq: convergence}
 \mb{E}\big(\mathop{\sup_{t \in [0, T]}|Y_{t}^{m} - \tilde{Y}_{t}|}\big) + |Z^{m} -\tilde{Z}|_{L^{2}(W)} +|U^{m} -\tilde{U}|_{L^{2}(\tilde{N}_{p})} \to 0, 
\end{equation}
and ($\tilde{Y}, \tilde{Z}, \tilde{U}$) solves the BSDE (\ref{eq: EDSRavecsauts}) given by ($f$, $B$).
 \end{proposition}
We first check that ($f^{m}$) constructed in Step 2 satisfies (1) in Lemma \ref{convmonotone} 
\[  \begin{array}{l}
\disp{| f^{m}(s, z^{m}, u^{m}) - f(s, z, u)| } \\
\\
  \quad \quad  \disp{ \leq \underbrace{|(f^{m} - f)(s, z^{m}, u^{m})|}_{= (I)} + \underbrace{|f(s, z^{m}, u^{m}) - f(s, z, u)|}_{= (II)}.}\\
 \end{array} \]
 Thanks to the continuity with respect to $z$ and $u$ of $f$, whose expression is (\ref{eq: generateurf}), (II) converges to zero. 
Furthermore, the boundedness of both sequences ($u^{m}$) and ($z^{m}$) ensures that, for $m$ large enough, $f$ and $f^{m}$ coincide: hence, (I) is equal to zero, which leads to the conclusion.\\  
\begin{flushright}
 $\square$
\end{flushright}

\textbf{Proof of Lemma \ref{convmonotone}} $\quad$ 
 We relegate to Appendix B the tedious proof of the strong convergence of ($Z^{m}$) and ($U^{m}$) in their respective Hilbert spaces and, assuming this, we prove the existence of a solution of the BSDE (\ref{eq: EDSRavecsauts}) with parameters ($f$, $B$).
 To identify ($\tilde{Y}$,$\tilde{Z}$, $\tilde{U}$) as a solution of the BSDE (\ref{eq: EDSRavecsauts}), we have to prove\\
 (i) $\disp{\int_{0}^{t}Z_{s}^{m}dW_{s}  \to  \int_{0}^{t}\tilde{Z}_{s}dW_{s},  \; \textrm{as} \; m \to \infty}$,\\
  (ii)$\disp{\int_{0}^{t}\int_{\mathbb{R}^{*}}U_{s}^{m}\tilde{N}_{p}(dx,ds) \to  \int_{0}^{t}\int_{\mathbb{R}^{*}}\tilde{U}_{s}\tilde{N}_{p}(dx,ds), \; \textrm{as} \; m \to \infty}$,\\
  (iii)  $\disp{\int_{0}^{t}f^{m}(s,Z_{s}^{m},U_{s}^{m})ds \to \int_{0}^{t}f(s,\tilde{Z}_{s},\tilde{U}_{s})ds, \; \textrm{as} \; m \to \infty}.$\\
Firstly, since $\tilde{Z}$ and $\tilde{U}$ are the weak limits of ($Z^{m}$) and ($U^{m}$),
 assertions (i) and (ii) correspond to the strong convergence, respectively, in $L^{2}(W)$ for (i) and, in $L^{2}(\tilde{N}_{p}(dx, ds))$ for (ii).\\
\indent Now, to prove (iii), we need to justify that the convergence holds true in $L^{1}(ds \otimes d\mb{P})$. This is argued as follows:\\
$\bullet \;$ on the one hand and from (i) and (ii), we have, eventually along a subsequence, the convergence in $ds \otimes d\mb{P}$-measure of ($Z^{m}$) and ($U^{m}$). 
From this remark and using assertion (1) in lemma \ref{convmonotone}, we deduce the convergence in $ds \otimes d\mb{P}$-measure of the sequence $(f^{m}(s, Z_{s}^{m}, U_{s}^{m}))$ to $f(s, \tilde{Z}_{s}, \tilde{U}_{s})$.\\  
$\bullet \;$ on the other hand,
we prove the uniform integrability of $(f^{m}(s, Z_{s}^{m}, U_{s}^{m}))$, along the subsequence where $(Z^{m})$ and $(U^{m})$ converge. Using now assumption ($H_{1}$) satisfied by each $f^{m}$ and Corollary \ref{equivalence} to argue that ($|U^{m}|_{\alpha}$) is uniformly bounded, this implies
$$\exists \; K^{1} \in L^{1}(ds \otimes d\mb{P}), \; |f^{m}(s, Z_{s}^{m}, U_{s}^{m})| \leq K_{s}^{1} + \alpha|Z_{s}^{m} - \tilde{Z}_{s}|^{2} + \alpha|\tilde{Z}_{s}|^{2} + C,$$
  with $K^{1}$ and $C$ depending only on the parameters $\alpha$ and $\theta$ appearing in ($H_{1}$).
The strong convergence in $L^{1}(ds \otimes d\mb{P})$ of $(|Z^{m} - \tilde{Z}|^{2})_{m}$ implying its uniform integrability, the result follows.\\
Using (i), (ii) and (iii), we obtain that ($\tilde{Y}, \tilde{Z}, \tilde{U}$) satisfies
 \begin{equation}\label{eq: eqoriginal}
\tilde{Y}_{t} = B + \int_{t}^{T}f(s,\tilde{Z}_{s},\tilde{U}_{s})ds - \int_{t}^{T}\tilde{Z}_{s}dW_{s} - \int_{t}^{T}\int_{\mathbb{R}^{*}}\tilde{U}_{s}\tilde{N}_{p}(dx,ds).
 \end{equation}
 To prove the convergence given by (\ref{eq: convergence}) in Lemma \ref{convmonotone}, we just take the supremum over $t$ and then the expectation in It\^o's formula applied to $\tilde{Y} - Y^{m}$ 
 \begin{tabbing} 
$\mb{E}\left(\disp{\sup_{t \in [0, T]}|\tilde{Y}_{t} - Y_{t}^{m}|}\right) $\= $\leq \; \disp{\mb{E}\left(\int_{0}^{T}|f(s,\tilde{Z}_{s},\tilde{U}_{s}) - f^{m}(s, Z_{s}^{m}, U_{s}^{m})|ds\right)} $\\
\\
\> $ \;+ \;\mb{E}\left(\disp{\sup_{t \in [0, T]}| \int_{t}^{T} (\tilde{Z}_{s} -Z_{s}^{m})dW_{s} |}\right) $\\ 
\\
\> $ \; + \;
\mb{E}\left(\disp{\sup_{t \in [0, T]}| \int_{t}^{T} ( \tilde{U}_{s} - U_{s}^{m})\tilde{N}_{p}(ds, dx)|}\right).$\\
\end{tabbing} 
Now, thanks to the convergence given by (iii) and Doob's inequality for the square integrable martingales, it follows that: $\mb{E}\left(\disp{\sup_{t}|\tilde{Y}_{t} - Y_{t}^{m}|}\right) \to 0$.
 \begin{flushright}
$\square$
           \end{flushright} 

\section{Application to finance: the compact case}
\subsection{The optimization problem}
The utility maximization problem considered is associated with the exponential utility function $U_{\alpha}$ with real parameter $\alpha$ $\big($ $U_{\alpha}(.):= -\exp(-\alpha \cdot)$ $\big)$. This problem consists in maximizing the expected value of the exponential utility of the portfolio (i.e the wealth at time $T$ minus a liability denoted by $B$). More precisely, we characterize in this section the expression of the value process $V^{B}$ which is defined at time $t$ by
 \begin{equation}\label{eq: pboptim} V_{t}^{B}(x) = \disp{\sup_{\pi \in \mc{A}_{t}}
  \mathbb{E}(U_{\alpha}(x + \int_{t}^{T}\pi_{s} \frac{dS_{s}}{S_{s-}} - B)| \mc{F}_{t}) }. 
  \end{equation}
$B$ stands for the contingent claim, which is assumed to be an $\mc{F}_{T}$-measurable random variable and $x$ is a constant standing for the  welath at time $t$. 

\subsection{Statement of the main result}

\begin{theorem}\label{expressionprocessV}
For any constant $x$, the expression of $V_{t}^{B}(x)$, as defined in (\ref{eq: pboptim}), is  
\begin{equation}\label{eq: fonctionV} \disp{ V_{t}^{B}(x) = -\exp(- \alpha(x - Y_{t})),} \end{equation}
where $Y_{t}$ is the first component of the solution ($Y, Z, U$) of the BSDE (\ref{eq: EDSRavecsauts}) given by the parameters ($f$, $B$), and whose generator $f$ is 
\begin{eqnarray}
 \nonumber \disp{ f(s,z,u) = \disp{\inf_{\pi \in
 \C}\left(\frac{\alpha}{2}|\pi\sigma_{s} - (z+ \frac{\theta}{\alpha})|^{2} + |u - \pi \beta_{s}|_{\alpha} \right) -
  \theta z - \frac{|\theta|^{2}}{2\alpha}}  }. \\
    \nonumber
\end{eqnarray}
Moreover, there exists an optimal strategy $\pi^{*}$ such that: $\pi^{*} \in \mc{A}_{t}$, and satisfying
 \begin{equation}\label{eq: optimstrat}
 \disp{\pi^{*}_{s} \in  \textrm{arg}\displaystyle{\min_{\pi \in \mc{C}}\left(\frac{\alpha}{2}|\pi\sigma_{s} - (Z_{s}+ \frac{\theta_{s}}{\alpha})|^{2} + |U_{s} - \pi \beta_{s}|_{\alpha}\right) }}.
\end{equation}
\end{theorem}

\subsection{Proof of theorem \ref{expressionprocessV}}
\paragraph*{The dynamic method}

\hspace{0.5cm}
As in \cite{ImkelleretHu}, the aim is to construct a family of processes ($R^{\pi}$)$_{\pi \in
\mathcal{A}_{t}}$, ($t$ being fixed), which is defined on $[t, T]$ and which satisfies
\[ \begin{array}{l}
(i) \; \disp{R_{t}^{\pi} = R_{t} \; \textrm{is a fixed}\;\mc{F}_{t}\; \textrm{-measurable random variable}\;, } \\
\\
(ii)\;  \disp{R_{T}^{\pi} = -\exp(-\alpha(X_{T}^{\pi} - B)),} \\
\\
(iii) \; \disp{R^{\pi} \; \textrm{is a supermartingale for each} \; \disp{\pi} \; \textrm{in}\; \mathcal{A}_{t},\; \textrm{and 
       there exists} \; \pi^{*},}\\ \disp{\pi^{*} \in \mathcal{A}_{t}}, \;   
\textrm{such that} \; \disp{R^{\pi^{*}}} \textrm{is a martingale.}   \\
\end{array}  \]
For this, we look for a process $R^{\pi}$ such that
$$\forall \; s \in [t, T], \quad R_{s}^{\pi} = U_{\alpha}(X_{s}^{\pi} - Y_{s}).$$ For sake of clarity, we use the notation $X^{\pi}$ instead of $X^{\pi, t, x} $ and we assume that there exists a triple ($Y$, $Z$, $U$)) solving a BSDE with jumps of the form (\ref{eq: EDSRavecsauts}),
with terminal condition $B$ and with a generator $f$ to be determined. 
 Referring to Theorem 5.1, Chapter 2 in \cite{Ikedawatanabe}), we first apply a generalized It\^o's formula to $R^{\pi}$, for any strategy $\pi$,
 \[ \begin{array}{l}
 \disp{   R_{s}^{\pi} - R_{t}^{\pi} =  - \alpha \int_{t}^{s}R_{u}^{\pi}(\pi_{u}\sigma_{u} - Z_{u})dW_{u} }   \\
\disp{ \quad + \int_{t}^{s}R_{u^{-}}^{\pi}\int_{\mathbb{R}^{*}}{(\exp(- \alpha(\pi_{u}\beta_{u} - U_{u})) - 1)\tilde{N}_{p}(du,dx)})} \\
 \disp{ \quad - \alpha \big(\int_{t}^{s}R_{u}^{\pi}(\pi_{u}b_{u} +
   f(u, Z_{u}, U_{u}))du\big)   + \frac{\alpha^{2}}{2}\int_{t}^{s}{R_{u}^{\pi}|\pi_{u}\sigma_{u} - Z_{u}|^{2}du} }  \\
  \disp{\quad + \int_{t}^{s} R_{u}^{\pi}\int_{\mathbb{R}^{*}}(\exp(- \alpha(\pi_{u}\beta_{u} - U_{u})) - 1 + \alpha(\pi_{u}\beta_{u} - U_{u}))\hat{N}_{p}(du,dx) }. \\
    \end{array} \]
 $R^{\pi}$ satisfies: $\disp{dZ = Z_{-}dM^{\pi} + ZdA^{\pi}}$, 
with $A^{\pi}$ such that
\[ \begin{array}{ll}
 dA_{u}^{\pi} := & \big(\frac{\alpha^{2}}{2}|\pi \sigma_{u} - Z_{u}|^{2} - \alpha(\pi b_{u} + f(u, Z_{u},U_{u}))\big)du \\
\\
   &\;  + \;\disp{ \int_{\mathbb{R}^{*}}g_{\alpha}(U_{u}(x))n(dx)du}, \\
   \end{array} \]
and with $\tilde{M}_{t, s}^{\pi}:= \frac{\mc{E}_{s}(M^{\pi})}{\mc{E}_{t}(M^{\pi})}$, where $\mc{E}(M^{\pi})$ denotes the Doleans-Dade exponential of the local martingale $M^{\pi}$
$$ M_{u}^{\pi} = \underbrace{(-\alpha (\pi \sigma - Z) \cdot W)_{u}}_{ \disp{= M_{u}^{1}}} + \underbrace{(\exp(- \alpha(\pi \beta -
  U)) - 1) \cdot \tilde{N}_{p})_{u}}_{ \disp{= M_{u}^{2}}} ,$$ 
with $M^{1}$ (resp. $M^{2}$) standing for the continuous part of $M^{\pi}$ (resp. the discontinuous part).
It follows that $R^{\pi} $ has the multiplicative form.
\begin{equation}\label{eq: multiplicativedecomp}
\disp{\forall \; s \ge t , \quad R_{s}^{\pi} = R_{t}^{\pi}\tilde{M}_{t, s}^{\pi}\exp(\disp{A_{s}^{\pi} -A_{t}^{\pi} }).}
\end{equation}
 Since: $(\exp(- \alpha(\pi \beta -
  U)) - 1) \ge -1, \; \mb{P}$-a.s., the Doleans-Dade exponential of $M^{2}$ is a positive local martingale and hence, a supermartingale. The supermartingale condition in (iii) holds true, provided, for all $\pi$, 
the process $\tilde{A}^{\pi}:= \exp(A^{\pi})$ is non decreasing: this entails  \\
    \begin{tabbing}
   $ \frac{\alpha^{2}}{2}|\pi \sigma_{u} - Z_{u}|^{2} $ \= $- \alpha(\pi b_{u} + f(u,Z_{u},U_{u})) $ \\
    \> $ + \disp{\int_{\mathbb{R}^{*}}{(\exp(- \alpha(\pi\beta_{u} - U_{u})) - 1 + \alpha(\pi\beta_{u}- U_{u}))n(dx)}
   \ge 0.}$ \\
   \end{tabbing}
This condition holds true, if we define $f$ as follows\\
\begin{eqnarray*}
 \disp{ f(s,z,u) = \displaystyle{\inf_{\pi \in
 \mathcal{C}}\left(\frac{\alpha}{2}|\pi\sigma_{s} - (z+ \frac{\theta}{\alpha})|^{2} + |u - \pi \beta_{s}|_ {\alpha}\right) -
  \theta_{s} z - \frac{|\theta_{s}|^{2}}{2\alpha}}  }. \\
   \nonumber
\end{eqnarray*}
Provided: $u \in L^{2} \cap L^{\infty}$, this assumption ensures that $|u -\pi\beta|_{\alpha} $ is finite for any $\pi \in \mathcal{A}_{t}$, thanks to the boundedness of $\beta$ and $\pi$ ($\pi$ taking its values in the compact set $\mathcal{C}$) and for any $z$, $z \in \mathbb{R}$, $f(s, z, u)$ is almost surely finite.\\

\paragraph*{Expression of the value function and optimal strategies}
To prove the supermartingale property of $R^{\pi}$ for any strategy $\pi$ ($\pi \in \mathcal{A}$), we rely on the results obtained in Section 3 on the BSDE with parameters ($f$, $B$). 
For this, we use the multiplicative form (\ref{eq: multiplicativedecomp}) of $R^{\pi}$ obtained in the previous paragraph.\\
\indent  Being a stochastic exponential of a martingale with jumps strictly larger than $-1$, $\mc{E}(M^{\pi})$ is a positive local martingale for any $\pi$, and consequently, there exists a sequence of stopping times $(\tau^{n})$ converging to $T$ such that $\tilde{M}_{. \wedge\tau^{n}}^{\pi}$ is a martingale. Since $e^{A^{\pi}}$ is non decreasing and $R_{t}$ is non positive, we can claim, on the one hand, that $R_{. \wedge\tau^{n}}^{\pi}$ satisfies
 \begin{equation}\label{eq: supermproperty} 
\forall \; s \le u, \;\forall A \in \mc{F}_{s}, \quad \mathbb{E}(R_{u \wedge \tau^{n}}^{\pi}\mathbf{1}_{A}) \leq \mathbb{E}(R_{s \wedge \tau^{n}}^{\pi}\mathbf{1}_{A}).
 \end{equation}
 On the other hand, since
$$ \;  R_{t}^{\pi} := - e^{-\alpha X_{t}^{\pi}}e^{\alpha Y_{t}} = - e^{-\alpha x}e^{\alpha Y_{t}}, $$
we use both the uniform integrability of $(e^{- \alpha X_{\tau}^{\pi}})$ (resulting from Lemma \ref{classequality}), where $\tau$ runs over the set of all stopping times and the boundedness of $Y$ to obtain the uniform integrability of ($R_{. \wedge\tau^{n}}^{\pi}$)$_{n}$. Hence, the passage to the limit as $n$ goes to $\infty$ in (\ref{eq: supermproperty}) is justified 
and it implies
$$\forall \; s, \; s \le u, \; \forall \;A \in \mathcal{F}_{s}, \quad \mathbb{E}(R_{u}^{\pi}\mathbf{1}_{A}) \leq \mathbb{E}(R_{s}^{\pi}\mathbf{1}_{A}), $$
 which entails the supermartingale property of $R^{\pi}$.\\
\indent To complete the proof, we justify the expression of $V^{B}(x)$ by proving the optimality of any strategy $\pi^{*}$ satisfying (\ref{eq: optimstrat}).
In fact, for this expression of $\pi^{*}$, we have: $A^{\pi^{*}} \equiv 0$ and hence, $R^{\pi^{*}} = R_{t}^{\pi^{*}}\tilde{M}^{\pi^{*}}$ is a true martingale ($\pi^{*}$ is in $\mathcal{A}_{t}$, thanks to Lemma 1). As a result, 
$$ \disp{\sup_{\pi \in \mc{A}_{t}} \mb{E}(R_{T}^{\pi})} = R_{t}^{\pi^{*}} = V_{t}^{B}(x). $$
Using that ($Y, Z, U$) is the unique solution of the BSDE given by ($f$, $B$), we obtain the expression (\ref{eq: fonctionV}) for the value function.\\ 
\begin{flushright}
$ \square$
\end{flushright}

\subsection{Characterization of optimal strategies}
The following lemma answers positively to both problems of existence and measurability of a strategy $\pi^{*}$ satisfying (\ref{eq: optimstrat}).
	 \begin{proposition}\label{propertyoptim}
	 Let $Z$ and $U$ be two predictable processes taking their values in $\mathbb{R}$ and $L^{2}(n(dx))$ and $\mathcal{C}$ a subset of $\mathbb{R}$.
	 \begin{enumerate}
\item[a.]  The process $F$ defined as below is again predictable 
$$\forall s \in \left[0 , T \right], \; F(s,Z_{s},U_{s}) =\disp{\inf_{\pi \in \mathcal{C} } \left(\frac{\alpha}{2}|\pi\sigma_{s} -
(Z_{s}+\frac{\theta}{\alpha})|^{2} + |U_{s} - \pi \beta_{s}|_{\alpha}\right).
  } $$
\item[b.] There exists a predictable version of $\pi^{*}$ which attains, for all $s$, $\omega$, the minimum taken over $\mathcal{C}$ of the sum of the two convex functionals\\   
	$$\frac{\alpha}{2}|\pi\sigma_{s} -
(Z_{s}+\frac{\theta}{\alpha})|^{2} \;\;  \textrm{and} \; \; |U_{s} - \pi \beta_{s}|_{\alpha}.$$
\end{enumerate}
\end{proposition}

\textbf{Proof of lemma \ref{propertyoptim}}
To prove assertion $a.$, we first introduce the sequence ($F^{n}$)$_{n}$ of 
predictable processes\\
\begin{tabbing}
$F^{n}$ \=$ :=F^{n}(s,Z_{s},U_{s}) $\\
\>$\;:= \disp{\inf_{\atop \{\pi^{n} \in \mb{Q}, \; \disp{dist(\pi^{n}, \mc{C}) \le \frac{1}{n}}\}} \left( \frac{\alpha}{2}|\pi^{n}\sigma_{s} -
(Z_{s}+\frac{\theta}{\alpha})|^{2} + |U_{s} - \pi^{n} \beta_{s}|_{\alpha}\right)}. $\\
\end{tabbing}
$F^{n}$ is a predictable process: in fact, on the one hand,
the infimum is taken over a countable subset of $\mb{Q}$. On the other hand, the functionals of the predictable processes $Z$ and $U$ are continuous. Besides and thanks to the fact that:
$F(s,Z_{s},U_{s})= \disp{ \lim_{n} F^{n}(s,Z_{s},U_{s})},$
 $F$ is itself predictable as a limit of such processes.\\
Now, to justify assertion $b.$, we argue the existence of a predictable selection by applying a measurable selection theorem (one reference is \cite{Aumann69}) to the set 
\begin{tabbing}
$  $ \= $ G := (t, \omega) \quad \textrm{s.t.}  $\\
\> $ \; \;\{ |f(t, \omega, \pi, Z_{t}(\omega), U_{t}(\omega)) - \disp{\inf_{\pi \in \mc{C}} f(t, \omega, \pi, Z_{t}(\omega), U_{t}(\omega))}| + (1 - \mathbf{1}_{\mc{C}}(\pi)) = 0\}.$
\end{tabbing}
The presence of the last term in the left-hand side ensures that the predictable choice of  $\pi :=(\pi(s, \omega))$ takes its values in $\mc{C}$.
 \begin{flushright}
$\square$
           \end{flushright}

\section{The non compact case}
To solve the problem (\ref{eq: pboptim}) and give the expression of the value function in the non compact case, we need to
solve the same BSDE with parameters ($f$, $B$) and with $f$ always given by (\ref{eq: generateurf}).
We again check that: $f:(z, u) \to f(s, z, u)$ is finitely valued ($\mb{P}$-a.s. and for all $s$): in fact, $f$ is a convex continuous functional of $z$ and $u$, which tends to $\infty$, as $|z|$ and $|u|_{L^{\infty}}$ goes to $\infty$ and hence, the infimum is attained. But, as soon as $\mc{C}$ is no more compact, the li-nearization procedure used to prove the uniqueness result cannot be applied: in this case, the BMO controls obtained in section 3.4 for the processes $\kappa$ and $\gamma$ do not hold anymore. This does not exclude a priori that there could be a unique solution.\\
\indent The method used as well as the justifications of the proof of the existence result being very similar as in the compact case, we just give here an outline of the procedure. The method consists once again in constructing a good approximation: for this, we introduce a sequence of BSDEs with parameters ($f^{m}$, $B$) such that:\\
\begin{itemize}
\item[1.] The sequence $(f^{m})_{m}$ is monotonic and it satisfies
$$f^{m}(s, z, u) \to  f(s, z, u), \; \mb{P}\textrm{-a.s. and for all} \; s,\; \textrm{as}\; m \to \infty .$$
\item[2.] All the BSDEs given by ($f^{m}$, $B$) satisfy ($H_{1}$) and ($H_{2}$): hence, thanks to Section 3, both the existence, uniqueness and comparison results hold for this sequence of BSDEs.
\end{itemize}
We then define ($f^{m}$)
\begin{eqnarray}\label{eq: generateurfm}
\nonumber \disp{ f^{m}(s,z,u) = \displaystyle{\inf_{\pi \in
 \mathcal{C}^{m}}\left(\frac{\alpha}{2}|\pi\sigma_{s} - (z+ \frac{\theta_{s}}{\alpha})|^{2} + |u - \pi \beta_{s}|_ {\alpha}\right)}}\\
  \disp{-\theta_{s} z - \frac{|\theta_{s}|^{2}}{2\alpha}}.\\
   \nonumber \end{eqnarray}
$\mc{C}^{m}$ being defined as follows: $\mc{C}^{m} := \mc{C} \cap [-m, m]$ corresponding to the intersection with the compact subset of $\mb{R}$.\\
\indent 
Thanks to the results obtained in Section 3, there exists a sequence ($Y^{m}, Z^{m}, U^{m}$) of solution of these BSDEs. Furthermore, since $(f^{m})$ is decreasing w.r.t. $m$,  it results from the comparison result, which is a direct byproduct of Theorem \ref{uniqueness}, that ($Y^{m}$) is decreasing ($\mb{P}$-a.s. and for all $s$).
Using now the results of Section 4, we get
$$ V_{m}^{B}(x) := U_{\alpha}(x - Y_{0}^{m}).$$
This corresponds to the optimization problem at time $t:=0$ and with $C^{m}$ as constraint set: ($Y_{0}^{m}$) being a decreasing sequence, ($V_{m}^{B}(x)$) is increasing, which implies 
$$\disp{\lim_{m} \nearrow \big(V_{m}^{B}(x)\big)} = \disp{\lim_{m} \nearrow \big(U_{\alpha}(x - Y_{0}^{m})\big)} = U_{\alpha}(x - Y_{0}). $$ 
\indent
The last step consists in identifying the limit ($Y, Z, U$) as a solution of the BSDE with parameters ($f$, $B$) and then conclude 
$$V^{B}(x) = U_{\alpha}(x - Y_{0}).$$ For this, we rewrite the proof of a ``monotone stability'' result analogous to lemma \ref{convmonotone}.
We just state the result, which is the key ingredient of the existence result.
\begin{proposition}\label{convmonotoneter} 
We denote by ($Y^{m}, Z^{m}, U^{m}$) the solutions of BSDEs with parameters ($f^{m}$, $B$) and we assume that ($f^{m}$) is such that:\\
 \textbullet $\;$ ($f^{m}$) satisfies the monotone convergence 
$$\disp{\lim_{m} \searrow f^{m}(s,Z_{s}^{m},U_{s}^{m}) }:= f(s,\tilde{Z}_{s},\tilde{U}_{s})\; \mb{P}\textrm{-a.s. and for all}\; s, $$
with $f$ continuous w.r.t. the variables $z$ and $u$.\\
  \textbullet $\;$ Each $f^{m}$ satisfies ($H_{1}$) (with parameters independent of $m$).\\
Then, ($Y^{m}$, $Z^{m}$, $U^{m}$) converges to ($\tilde{Y}, \tilde{Z}, \tilde{U}$) in the following sense
\begin{equation}\label{eq: Convergence}
\mb{E}(\disp{\sup_{t \in [0, T]}|Y_{t}^{m} - \tilde{Y}_{t}|}) + |Z^{m} -\tilde{Z}|_{L^{2}(W)} +|U^{m} -\tilde{U}|_{L^{2}(\tilde{N}_{p})} \to 0. \end{equation}
 $(Y^{m})$ is decreasing and the triple ($\tilde{Y}, \tilde{Z}, \tilde{U}$) is in $S^{\infty} \times L^{2}(W) \times L^{2}(\tilde{N}_{p})$ and solves the BSDE with parameters ($f$, $B$).
 \end{proposition}
 For the rest of the proof, we refer the reader to Paragraph 3.5.3 and, in particular, to Appendix B for the strong convergence given by (\ref{eq: Convergence}). The proof of this last convergence is the key point: we refrain from doing this since the proof given in Appendix B can be rewritten identically (if we replace the increasing sequence of generators by a decreasing one) and, in particular, the estimates of lemma \ref{estim2} hold true for the sequences ($Y^{m}$), ($Z^{m}$) and ($ U^{m}$).\\
  \begin{flushright}
$\square$
           \end{flushright}
\section{Conclusion}
In that paper, we have solved the utility maximization problem with portfolio constraints in a discontinuous setting and by means of BSDEs. Our setting is very similar as the one in \cite{Becherer}, but, contrary to the author of the aforementionned paper, the presence of constraints in the model entails that the driver has quadratic growth w.r.t. its variable $z$. Hence, our main achievement consists in obtaining existence and uniqueness results for such quadratic BSDEs with jumps. This theoretical study and the use of the dynamic principle allow to give the expression of the value function for any time $t$ in terms of the solution of the BSDE and characterize optimal strategies for the problem.\\
Due to the restrictions, some questions remain for further investigations: in particular, the case when the Levy measure is infinite is unsolved (this assumption is already given in \cite{Becherer}). Here, we also restrict our study to the exponential utility function and we mention here that, in the power utility case, the utility maximization problem with non zero liability is an open problem.
\\

\subsection{Omitted proofs}
\subsection{Appendix A: proof of lemma \ref{estimcaslip}}
For convenience for the reader, we give a detailed outline of the proof which is adapted from Proposition 2.2 in \cite{BriandetCoquet} to the discontinuous setting.
Contrary to lemma \ref{estim2}, where the first component of the solution is supposed to be in $S^{\infty}$, in this proposition, $Y^{n}$ is only assumed to be in $S^{2}$. We first write It\^o's formula for $(e^{\Gamma t} |Y_{t}^{n}|^{2})$, 
 $\Gamma$ being a non negative constant which is explicited during the proof.
  \begin{equation}\label{eq: formuleIto1}
d(e^{\Gamma t} |Y_{t}^{n}|^{2}) = \Gamma e^{\Gamma t} |Y_{t}^{n}|^{2}dt + e^{\Gamma t}\big(2Y_{t}^{n} dY_{t}^{n} + d\lan{Y^{n}}\ran_{t} \big), 
  \end{equation}
with
 \begin{tabbing}
 $ 2Y_{t}^{n} dY_{t}^{n} + d\lan{Y^{n}}\ran_{t} $\=$:= -2Y_{t}^{n}g^{n}(t, Z_{t}^{n}, U_{t}^{n})dt  + \left(|Z_{t}^{n}|^{2} + \disp{\int_{\mb{R}\setminus\{0\}}|U_{t}^{n}(x)|^{2}n(dx)}\right)dt$\\
 \\
$ $ \> $\; \quad \;+  \;2Y_{t}^{n} \big(Z_{t}^{n}dW_{t} + \disp{\int_{\mb{R}\setminus\{0\}}U_{t}^{n}\tilde{N}_{p}(dt, dx)} \big). $ \\
 \end{tabbing}
Since ($Y^{n}, Z^{n}, U^{n}$) is in $S^{2} \times L^{2}(W) \times L^{2}(\tilde{N}_{p})$, the process $K$, such that\\
\begin{equation}\label{eq: formeK}
\forall \; s \in  [0,T], \quad K_{s}:= \disp{\int_{0}^{s}2 e^{\Gamma u}Y_{u}^{n} \big(Z_{u}^{n}dW_{u} + \disp{\int_{\mb{R}\setminus\{0\}}U_{u}^{n}\tilde{N}_{p}(du, dx)}\big)}, 
\end{equation}
 is a true martingale.
 We fix $t$ ($t \in [0,  T]$)
and we rewrite the formula (\ref{eq: formuleIto1}) in the integrated form between $s$ and $T$, for any $s$, $t \le s \le T$ , 
  \[ \begin{array}{l}
 e^{\Gamma s} |Y_{s}^{n}|^{2} -e^{\Gamma T} |Y_{T}^{n}|^{2} \\
 =\;\disp{\int_{s}^{T}e^{\Gamma u} Y_{u}^{n}\big(-\Gamma Y_{u}^{n} + \;2g^{n}(u, Z_{u}^{n}, U_{u}^{n})\big)du }\\

 \quad \;- \disp{\int_{s}^{T}e^{\Gamma u}\left(|Z_{u}^{n}|^{2} +\disp{\int_{\mb{R}\setminus\{0\}}|U_{u}^{n}(x)|^{2}n(dx)}\right)du -\big(K_{T} - K_{s} \big). } \\
  \end{array} \]
 We then rely on the Lipschitz property of the generator $g^{n}$  
 \begin{eqnarray}\label{eq: lipproperty}
 \nonumber \\ 
 2|Y_{u}^{n}||g^{n}(u, Z_{u}^{n}, U_{u}^{n} )| \leq 2|Y_{u}^{n}||g^{n}(u, 0, 0)| + 2L_{n}\big( |Y_{u}^{n}||Z_{u}^{n}| +    
|Y_{u}^{n}||U_{u}^{n}|_{L^{2}(n)}\big),\\
\nonumber
\\ 
\nonumber\end{eqnarray}
with: $|U_{u}^{n}|_{L^{2}(n)}:= \disp{\left( \int_{\mb{R}\setminus\{0\}}|U_{u}^{n}(x)|^{2}n(dx)\right)^{\frac{1}{2}}},$
 and using that: $|2L_{n}ab| \le ( 2(L_{n})^{2}a^{2} + \frac{1}{2}b^{2})$, we obtain
 \begin{equation}\label{eq: controle} 2L_{n}\big(|Y_{u}^{n}||Z_{u}^{n} | + |Y_{u}^{n}||U_{u}^{n}|_{L^{2}(n)}\big) \le 4(L_{n})^{2}|Y_{u}^{n}|^{2} + \frac{1}{2}\big(|Z_{u}^{n}|^{2} + |U_{u}^{n}|_{L^{2}(n)}^{2}\big).  \end{equation}
Combining both (\ref{eq: lipproperty}) and (\ref{eq: controle}), setting: $\Gamma = 4(L_{n})^{2}$ and taking the expectation w.r.t $\mc{F}_{t}$ in It\^o's formula applied to $ e^{\Gamma s} |Y_{s}^{n}|^{2}$ between $t$ and $T$, we get
 \begin{tabbing}
$ $ \=$\disp{ e^{\Gamma t} |Y_{t}^{n}|^{2}  \le} \disp{\mb{E}\left( e^{\Gamma T} |Y_{T}^{n}|^{2} |\mc{F}_{t}\right) }$\\
\\
$ $\> $ \; +\; \disp{\mb{E}\left( \int_{t}^{T}e^{\Gamma u}\left(2|Y_{u}^{n}||g^{n}(u, 0, 0)|  + \frac{1}{2}(|Z_{u}^{n}|^{2})+ \frac{1}{2}\disp{\int_{\mb{R}\setminus\{0\}} |U_{u}^{n}(x)|^{2} n(dx)} \right) du |\mc{F}_{t}\right)} $ \\
\\
$ $\> $ \;  - \disp{\mb{E}\left(\int_{t}^{T}e^{\Gamma u}\big(|Z_{u}^{n}|^{2}du +|U_{u}^{n}|_{L^{2}(n)}^{2} du \big)|\mc{F}_{t}\right)}.$ \\
\end{tabbing}
This leads to\\
\newline
$\disp{\mb{E}\left(\int_{t}^{T}e^{\Gamma u}\big(|Z_{u}^{n}|^{2} + |U_{u}^{n}|_{L^{2}(n)}^{2}\big) du  \right)| \mc{F}_{t}\big)}$\\
\begin{equation}\label{eq: Ito2} \quad \quad \leq \disp{2\left(\mb{E}\big(e^{\Gamma T}|Y_{T}^{n}|^{2} + 2\int_{t}^{T} e^{\Gamma u}|Y_{u}^{n}||g^{n}(u, 0, 0)| du| \mc{F}_{t}\big)\right).}
\end{equation}
We consider again the integrated form of (\ref{eq: formuleIto1}) between $s$ and $T$ and then, we take the supremum over $s$ ($s \in [t, T]$)
\[ \begin{array}{l}
 \disp{\sup_{ t\le s \le T}e^{\Gamma s} |Y_{s}^{n}|^{2}} \le    \disp{e^{\Gamma T} |Y_{T}^{n}|^{2}}
\\ \quad  \;+ \;\disp{2\int_{t}^{T}e^{\Gamma u}|Y_{u}^{n}||g^{n}(u, 0, 0)| du} 
\; + \disp{\sup_{ t\le s\le T}|K_{T} - K_{s}|} .\\
\end{array} \]
We now apply the Burkholder-Davis-Gundy inequality to the supremum of the square integrable martingale $K$ and the relation: $ Cab \le \frac{C^{2}}{2}a^{2} +\frac{1}{2}b^{2}$, to deduce
the existence of a constant $C$ such that
\begin{tabbing}
$ $ \= $\mb{E}\left(\disp{\sup_{ t\le s \le T}e^{\Gamma s} |Y_{s}^{n}|^{2}} |\mc{F}_{t} \right) \le \mb{E}\left( e^{\Gamma T} |Y_{T}^{n}|^{2} + 2 \disp{\int_{t}^{T}e^{\Gamma u}|Y_{u}^{n}||g^{n}(u, 0, 0)| du|\mc{F}_{t} }\right)$\\
\\
$ $\>  $ \quad \quad + \frac{C^{2}}{2}\disp{\mb{E}\big(\int_{t}^{T}e^{\Gamma u}\big( |Z_{u}^{n}|^{2} + |U_{u}^{n}|_{L^{2}(n)}^{2} \big)du |\mc{F}_{t}\big)  } + \frac{1}{2}\disp{\mb{E}\big(\disp{\sup_{ t\le s \le T}e^{\Gamma s} |Y_{s}^{n}|^{2}} |\mc{F}_{t}\big)  }.$\\
\end{tabbing}
From now, this constant $C$ is generic and may vary from line to line. 
Combining this previous inequality with (\ref{eq: Ito2}), we deduce 
\begin{eqnarray*}
 \mb{E}\left(\disp{\sup_{ t\le s\le T}e^{\Gamma s} |Y_{s}^{n}|^{2}} +\disp{ \int_{t}^{T}e^{\Gamma u}\big( |Z_{u}^{n}|^{2} +  |U_{u}^{n}|_{L^{2}(n)}^{2}\big)du}|\mc{F}_{t} \right) \\
 \quad \quad \quad \le C \mb{E}\left(e^{\Gamma T} |Y_{T}^{n}|^{2} +  \disp{\int_{t}^{T}e^{\Gamma u}|Y_{u}^{n}||g^{n}(u, 0, 0)| du|\mc{F}_{t} }\right).\\
\end{eqnarray*}
To obtain the desired relation, we check 
\begin{tabbing}
$ $\=  $C\mb{E}\left(\disp{\int_{t}^{T}e^{\Gamma u}|Y_{u}^{n}||g^{n}(u, 0, 0)|du|\mc{F}_{t}} \right)$ \\   \\
\> $ \quad \;\le \; \frac{1}{2}\mb{E}\left(\disp{\sup_{ t\le u\le T}e^{\Gamma u} |Y_{u}^{n}|^{2}|\mc{F}_{t}} \right) + \frac{C^{2}}{2}\mb{E}\left(\disp{\big(\int_{t}^{T}e^{\frac{\Gamma}{2} u}|g^{n}(u, 0, 0)|du\big)^{2}|\mc{F}_{t}} \right).$\\

\end{tabbing}
Noting that: $|Y_{t}^{n}|^{2}  \le \mb{E}\left(\disp{\sup_{ t\le s\le T}e^{\Gamma s} |Y_{s}^{n}|^{2}} | \mc{F}_{t}\right), $ the result follows.\\ 

\subsection{Appendix B: End of the proof of Lemma \ref{convmonotone}}
\hspace{0.5cm} To prove the essential result stated in Lemma \ref{convmonotone}, i.e. the strong convergence of ($Z^{m}$) and ($U^{m}$), we adapt the method already used in \cite{mkobylanski}. To this end,
we consider the function $\Phi_{K}$ 
\[ \Phi_{K}(x) = \frac{e^{2Kx}- 2K x - 1 }{2K} \; \; \big(= g_{2K}(x) \big), \]
which is twicely continuously differentiable and satisfies
\[ \left\{  \begin{array}{l} 
\Phi_{K}(0) = 0,\; \textrm{and } \;\Phi_{K},\;\Phi_{K}^{''}  \ge 0,\\
 \Phi_{K}^{'}(x) \ge 0, \; \textrm{if}\; x \ge 0, \\
\Phi_{K}^{''} - 2K\Phi_{K}^{'} = 2K. \\
\end{array} \right. \] 
 These properties imply that, taking $m, p$ such that $m \geq p \geq M$, with $M$ given as in Step 1 in the proof of existence, we have: $\Phi_{K}^{'}(Y^{m}_{s} - Y^{p}_{s}) \ge$ 0, for all $s$ and $\mathbb{P}$-a.s.
In the sequel, we note $Y^{m, p}$ instead of $Y^{m} - Y^{p}$ (the same holds for $Z^{m, p}$ and $U^{m, p}$).
We now apply It\^o's formula to $\Phi_{K}(Y^{m, p})$ and we take the expectation between 0 and $T$ \\  
$  \disp{\mathbb{E}{\phi_{K}(Y_{0}^{m, p})}}  =  \disp{ \mathbb{E}\int_{0}^{T}(e^{2K Y_{s}^{m, p}}- 1)(f^{m}(s,Z_{s}^{m},U_{s}^{m})- f^{p}(s,Z_{s}^{p},U_{s}^{p}))ds}$
\begin{equation}\label{eq: formuleito}  \quad \;\disp{- \mathbb{E}\int_{0}^{T}K e^{2K Y_{s}^{m, p}}|Z_{s}^{m, p}|^{2}ds
 - \mathbb{E}\int_{0}^{T}e^{2K Y_{s}^{m, p}}\int_{\mathbb{R}^{*}} g_{2K}(U_{s}^{m, p}(x)) n(dx)ds}.
\end{equation}
 
($\Phi_{K}(Y^{m, p}) $ being a uniformly bounded process, the expectation of the martingale part vanishes).
To give an upper bound of  $$F^{m, p} = f^{m}(s, Z_{s}^{m}, U_{s}^{m}) - f^{p}(s, Z_{s}^{p}, U_{s}^{p}),$$ we rely on ($H_{1}$) and on the result (\ref{eq: relationeq}) given in Corollary \ref{equivalence} to claim 
\begin{eqnarray*} \exists \; C, \; f^{m}(s, Z_{s}^{m}, U_{s}^{m}) \leq \frac{\alpha}{2}|Z_{s}^{m}|^{2} + |U_{s}^{m}|_{\alpha} \leq \frac{\alpha}{2}|Z_{s}^{m}|^{2} + C |U_{s}^{m}|_{L^{2}(n)}^{2}.
\end{eqnarray*}
Thanks to assertion (i)(b) in Corollary \ref{equivalence}, there exists a constant always denoted by $C$ and depending only on $|Y^{m}|_{S^{\infty}}$, and on the parameters $\alpha$ and $n(\mb{R}\setminus\{0\})$) such that 
\begin{equation}\label{eq: controlfm} f^{m}(s, Z_{s}^{m}, U_{s}^{m})  \leq \frac{\alpha}{2}|Z_{s}^{m}|^{2} + C.
\end{equation}
 We also use that:
$|\theta_{s} Z_{s}| \leq \frac{1}{\alpha}|\theta_{s}|^{2}  + \frac{\alpha}{4}|Z_{s}|^{2},$
to obtain the existence of $\hat{C} \in L^{1}(ds \otimes d\mathbb{P})$ such that\\
\begin{equation}\label{eq: controlfp} - f^{p}(s, Z_{s}^{p}, U_{s}^{p}) \le \hat{C}_{s} + \frac{\alpha}{4}|Z_{s}^{p}|^{2}, \end{equation}
with $\hat{C}$ given by $\hat{C}_{s} = \frac{3}{2\alpha}|\theta_{s}|^{2} + \frac{|\theta_{s}|^{2}}{2\alpha}$.
 Using the convexity of: $z \to |z|^{2}$, we obtain\\
\[   \begin{array}{ll} \frac{\alpha}{2}|Z_{s}^{m}|^{2}  &
 \leq  \frac{\alpha}{2}\left(\frac{1}{3}|3 Z_{s}^{m, p}|^{2} +\frac{1}{3}|3(Z_{s}^{p}- \tilde{Z}_{s})|^{2} + \frac{1}{3}|3\tilde{Z}_{s}|^{2} \right)   \\
\\
 & \leq  \frac{3\alpha}{2}(|Z_{s}^{m, p}|^{2} + |Z_{s}^{p}- \tilde{Z}_{s}|^{2} + |\tilde{Z}_{s}|^{2}), \\
\end{array} \]
and, similarly,
\[   \begin{array}{l}
\disp{\frac{\alpha}{4}|Z_{s}^{p}|^{2} \leq \frac{\alpha}{2}\big(|Z_{s}^{p} - \tilde{Z}|^{2} + |\tilde{Z_{s}}|^{2}\big)}.\\
\end{array} \]
These estimates lead to 
$$F^{m,\;p} \le \hat{C}_{s} + 2\alpha \big( |Z_{s}^{m, p}|^{2} + |Z_{s}^{p} - \tilde{Z}|^{2} + |\tilde{Z_{s}}|^{2}\big) + C.$$
 
We set: $K = 4 \alpha$ and, transferring into the left-hand side the terms containing $|Z^{m, p}|^{2}$ or $|U^{m, p}|_{8\alpha}$, we rewrite the equation (\ref{eq: formuleito})
\begin{tabbing}
 \= $ \quad     \disp{\mathbb{E}{\phi_{K}(Y_{0}^{m, p})}}   + \disp{\mb{E}\int_{0}^{T}e^{8\alpha Y_{s}^{m, p}} |U_{s}^{m, p}|_{8 \alpha} ds}$ \\
$ $\\
 \> $\quad + \disp{\mb{E}\int_{0}^{T}2 \alpha e^{8\alpha Y_{s}^{m, p}}|Z_{s}^{m, p}|^{2}ds + \mb{E}\int_{0}^{T} 2 \alpha|Z_{s}^{m, p}|^{2}ds } $\\
\\
 \> $\quad \quad \quad \leq \disp{\mb{E}\int_{0}^{T}(e^{8\alpha Y_{s}^{m, p}} - 1)(\hat{C}_{s} + C + 2\alpha (|Z_{s}^{p} - \tilde{Z}_{s}|^{2} + |\tilde{Z_{s}}|^{2}))ds}. $\\
\end{tabbing}
We rely once again on the result given in Corollary \ref{equivalence} to claim the existence of a constant $C$
$$ \frac{1}{C}\disp{\mb{E}\int_{0}^{T}e^{8\alpha Y_{s}^{m, p}} |U_{s}^{m, p}|_{L^{2}(n)}^{2} ds} \le \disp{\mb{E}\int_{0}^{T}e^{8\alpha Y_{s}^{m, p}} |U_{s}^{m, p}|_{8 \alpha} ds}.$$
This entails, taking the limit inf as $m$ goes to $\infty$ and $p$ being fixed 

\begin{tabbing}
 \= $ \quad     \disp{\mathbb{E}{\phi_{K}( \tilde{Y}_{0} - Y_{0}^{ p})}} + \;\frac{1}{C}\disp{\lim{ \inf_{m}\mb{E}\int_{0}^{T}e^{8\alpha Y_{s}^{m, p}} |U_{s}^{m, p}|_{L^{2}(n)}^{2} ds}} $ \\
$ $\\
 \> $\quad + \disp{ \lim{ \inf_{m} \left( \mb{E}\int_{0}^{T} 2 \alpha e^{8\alpha Y_{s}^{m, p}}|Z_{s}^{m, p}|^{2}ds + \mb{E}\int_{0}^{T} 2 \alpha|Z_{s}^{m, p}|^{2}ds \right)}} $\\
\\
 \> $\quad \quad \quad \leq \disp{\mb{E}\int_{0}^{T}(e^{8\alpha (\tilde{Y}_{s} - Y_{s}^{ p})} - 1)(\hat{C}_{s} + C + 2\alpha (|Z_{s}^{p} - \tilde{Z}_{s}|^{2} + |\tilde{Z}_{s}|^{2}))ds}. $\\
\end{tabbing}
where the limit, as $m$ goes to $\infty$ ($p$ being fixed), in the right-hand side of this inequality results from Lebesgue's convergence theorem: in fact, we have that $(Y_{s}^{m})$ converges to $\tilde{Y}_{s}$, $\mb{P}$-a.s. and for all $s$, and: $\hat{C}$, $|Z^{p} - \tilde{Z}|^{2}$ and $|\tilde{Z}|^{2}$ are in $L^{1}(ds \otimes d\mb{P})$.
Before justifying the passage to the limit as $p$ goes to $\infty$, we use that $\tilde{Z}$ and $\tilde{U}$ are the respective weak limits of ($Z^{m}$) and ($U^{m}$) to obtain
$$\disp{\mb{E}\int_{0}^{T}\left(2 \alpha (e^{8\alpha (\tilde{Y}_{s} - Y_{s}^{p})} + 1)\right) |\tilde{Z}_{s} - Z_{s}^{p} |^{2} ds} \leq  \disp{ \lim{ \inf_{m} \mb{E}\int_{0}^{T} 2 \alpha \left(e^{8\alpha Y_{s}^{m, p}} +1 \right)|Z_{s}^{m, p}|^{2}ds  }},$$
 and
$$ \disp{\mb{E}\int_{0}^{T}|\tilde{U}_{s} - U_{s}^{p}|_{L^{2}(n)}^{2} ds} \leq \disp{\lim{ \inf_{m}\mb{E}\int_{0}^{T}e^{8\alpha (Y_{s}^{m, p})} |U_{s}^{m, p}|_{L^{2}(n)}^{2}} ds},$$
 (since: $e^{8\alpha Y_{s}^{m, p}} \geq 1$). Transferring now in the left-hand side of the last inequality of page 30 the unique term containing $|\tilde{Z} - Z^{p}|^{2}$, we get
\begin{tabbing}
\= $ \disp{ \mb{E}\int_{0}^{T}4 \alpha|\tilde{Z}_{s} - Z_{s}^{p}|^{2}ds
+ \frac{1}{C}\mb{E}\int_{0}^{T}|\tilde{U}_{s} - U_{s}^{p}|_{L^{2}(n)}^{2}ds }$\\
\\ 
\> $\quad \quad \disp{\leq \mb{E}\int_{0}^{T} (e^{8\alpha (\tilde{Y}_{s} - Y_{s}^{ p})} - 1)(\hat{C}_{s} + C + 2\alpha |\tilde{Z}_{s}|^{2}) ds} 
$ \\
\end{tabbing}
and hence, the desired convergence result follows: in fact, the limit in the right-hand side, as $p$ goes to $\infty$, is equal to zero, thanks to the Lebesgue's convergence theorem.\\

\begin{flushright}
$\square$
           \end{flushright} 


\begin{thebibliography}{1}
\bibitem{Aumann69}
Aumann, R.,
\newblock{ \em Measurable utility and the measurable choice theorem}
\newblock{ \em Editions du Centre Nat. Recherche Sci.}: 15--26, 1969.

\bibitem{BarlesBuck}
Barles, G., Buckdahn, R. and Pardoux, E.,
 \newblock{ \em Backward stochastic differential equations and
              integral-partial differential equations},
\newblock{ \em Stoch. Stoch. Rep..}, \textbf{60} : 57--83, 1997.
\bibitem{Becherer}
Becherer, D.,
\newblock{\em Bounded solutions to Backward SDE's with jumps for utility optimization and indifference hedging},
\newblock{ \em Ann. Appl. Probab.}, \textbf{16}(4) : 2027--2054, 2006.
\bibitem{BriandetCoquet}
Briand, P., Coquet, F., Hu, Y., M{\'e}min, J. and Peng, S.,
 \newblock {\em A converse comparison theorem for {BSDE}s and related
              properties of {$g$}-expectation},
\newblock {\em    Electron. Comm. Probab.}, \textbf{5} : 101--117, 2000.


\bibitem{BriandetHu}
 Briand, P. and Hu, Y.,
 \newblock {\em BSDE with quadratic growth and unbounded terminal value},
 \newblock {\em Probab. Theory Related Fields}, \textbf{136}(4) : 604--618, 2006.



\bibitem{DelbetSchach}
Delbaen, F. and Schachermayer, W.,
\newblock The mathematics of arbitrage,
\newblock {\em Springer Finance}, Springer-Verlag, Berlin, 2006.


 \bibitem{Dellacheriemeyer}
Dellacherie, C. and Meyer, P.-A.
\newblock {\em  Probabilit\'es et Potentiel. Th\'eorie des martingales. {C}hapitres {V} \`a {VIII}}, Hermann, 1980.


\bibitem{ElKaretPeng}
El Karoui, N., Peng, S. and Quenez, M.C.,
\newblock {\em Backward stochastic differential equations in finance},
\newblock{\em Math. Finance}, \textbf{7}(1) : 1--71, 1997.



 \bibitem{Fritelli}
Biagini, S. and Frittelli, M.,
  \newblock {\em On the super replication price of unbounded claims},
  \newblock {\em  Ann. Appl. Probab.}, \textbf{14}(4) : 1970--1991, 2004.

\bibitem{Follmer}
F{\"o}llmer, H. and Schied, A.,
    \newblock {\em An introduction in discrete time stochastic finance},
  \newblock {\em  de Gruyter}, Berlin, 2002.
   
\bibitem{ImkelleretHu}
Hu, Y., Imkeller, P. and M\"uller, M.,
\newblock{\em Utility maximization in incomplete markets},
\newblock{\em Ann. Appl. Probab.}, \textbf{15}(3) : 1691--1712, 2005.

\bibitem{Ikedawatanabe}
 Ikeda, N. and Watanabe, S.,
\newblock {\em Stochastic differential equations and diffusion processes},
\newblock {\em North-Holland Publishing Co.}, Amsterdam, 1989.

\bibitem{mkobylanski}
Kobylanski, M.,
\newblock {\em Backward stochastic differential equations and partial
              differential equations with quadratic growth},
\newblock {\em Ann. Probab.}, \textbf{28}(2) : 558--602, 2000.


 \bibitem{Royer2}
Royer, M.,
  \newblock{\em Backward stochastic differential equations with jumps and related non-linear expectations},
 \newblock {\em Stochastic Process. Appl.}, \textbf{116}(10) : 1358--1376, 2006.


 \bibitem{Kazamaki} 
Kazamaki, N.,  
\newblock {\em  A sufficient condition for the uniform integrability of
              exponential martingales},
\newblock {\em Math. Rep. Toyama Univ.}, \textbf{2} : 1--11, 1979.
  
\bibitem{Pardoux}
Pardoux, {\'E}.,
 \newblock {\em Generalized discontinuous backward stochastic differential
              equations},
\newblock {\em Backward stochastic differential equations},
   Pitman Res. Notes Math. Ser.,
    \textbf{364} : 207--219, 1997.

\bibitem{Protter}
Protter, P., 
\newblock {\em Stochastic integration and differential equations}, Springer, Berlin, 2004.

\bibitem{Schach1}
Schachermayer, W.,
\newblock {\em Utility maximization in incomplete markets},
\newblock {\em Lecture Notes in Math.}, \textbf{1856} : 255--293, Springer, Berlin, 2004.

 \end{thebibliography}
\end{document}